\documentclass{gtart}

\def\ifplaintex{\expandafter\ifx\csname documentclass\endcsname\relax}

\def\gtp{{\mathsurround=0pt\it $\cal G\mskip-2mu$eometry \&\ 
$\cal T\!\!$opology $\cal P\!$ublications}}  

\def\recd{{\small Received:\qua\receiveddate\ifx\reviseddate\relax
\else\qquad Revised:\qua\reviseddate\fi\par}} 

\def\lognumber#1{\def\thelognumber{#1}}
\def\volumenumber#1{\def\thevolumenumber{#1}}
\def\volumeyear#1{\def\thevolumeyear{#1}}
\def\papernumber#1{\def\thepapernumber{#1}}
\def\pagenumbers#1#2{\def\startpage{#1}\def\finishpage{#2}}
\def\published#1{\def\publishdate{#1}}

\def\received#1{\def\receiveddate{#1}}
\def\revised#1{\def\reviseddate{#1}}
\def\accepted#1{\def\accepteddate{#1}}
\def\asciititle#1{\def\theasciititle{#1}}
\def\covertitle#1{\def\thecovertitle{#1}}

\long\def\asciiabstract#1{\long\def\theasciiabstract{#1}}

\let\\\par\let\thelognumber\relax\let\thevolumenumber\relax
\let\thepapernumber\relax\let\thevolumeyear\relax\let\startpage\relax
\let\finishpage\relax\let\publishdate\relax\let\receiveddate\relax
\let\reviseddate\relax\let\accepteddate\relax\let\theasciititle\relax
\let\thecovertitle\relax\let\theasciiauthors\relax
\let\theasciiabstract\relax

\let\theasciiemail\relax

\ifplaintex
\font\logobig=cmssbx10 scaled 3836
\font\logomed=cmssbx10 scaled 2557
\else
\font\logobig=cmssbx10 scaled 4200
\font\logomed=cmssbx10 scaled 2800
\fi

\long\def\makeagttitle{   
\count0=\startpage
\agt\hfill      
\hbox to 45truept{\vbox to 0pt{\vglue -13truept{\logomed A\kern -.37em{\logobig 
T}\kern -.38em G}\vss}\hss}
\break
{\small Volume \thevolumenumber\ (\thevolumeyear)
\startpage--\finishpage\nl
Published: \publishdate}

\vglue .25truein

{\parskip=0pt\leftskip 0pt plus
1fil\def\\{\par\smallskip}{\Large\bf\thetitle}\par\medskip} \vglue
0.05truein

{\parskip=0pt\leftskip 0pt plus 1fil\def\\{\par}{\sc\theauthors}
\par\medskip}%
 
\vglue 0.03truein 

{\small\leftskip 25truept\rightskip 25truept{\bf Abstract}\stdspace\theabstract

{\bf AMS Classification}\stdspace\theprimaryclass
\ifx\thesecondaryclass\relax\else; \thesecondaryclass\fi\par
{\bf Keywords}\stdspace \thekeywords\par}\vglue 7truept

}   

\ifplaintex
\font\phead=cmsl9 scaled 950
\font\pnum=cmbx10 scaled 913
\font\pfoot=cmsl9 scaled 950
\headline{\vbox to 0pt{\vskip -4.5mm\line{\small\phead\ifnum
\count0=\startpage ISSN 1472-2739 (on-line) 1472-2747 (printed)
\hfill {\pnum\folio}\else\ifodd\count0\def\\{ }%
\ifx\theshorttitle\relax\thetitle\else\theshorttitle\fi\hfill{\pnum\folio}
\else\def\\{ and }{\pnum\folio}\hfill\ifx\theshortauthors\relax\theauthors
\else\theshortauthors\fi\fi\fi}\vss}}
\footline{\vbox to 0pt{\vglue 0mm\line{\small\pfoot\ifnum\count0=\startpage
\copyright\ \gtp\hfill\else
\agt, Volume \thevolumenumber\ (\thevolumeyear)\hfill\fi}\vss}}
\else
\font=cmsl9 scaled 1050
\font\lnum=cmbx10 
\font=cmsl9 scaled 1050
\makeatletter
\def\@oddhead{{\small\lhead\ifnum\count0=\startpage ISSN 1472-2739 
(on-line) 1472-2747 (printed)\hfill {\lnum\number\count0}\else\ifodd\count0
\def\\{ }\ifx\theshorttitle\relax \thetitle \else\theshorttitle\fi\hfill
{\lnum\number\count0}\else\def\\{ and }{\lnum\number\count0}
\hfill\ifx\theshortauthors\relax 
\theauthors\else\theshortauthors\fi\fi\fi}}\def\@evenhead{\@oddhead}
\def\@oddfoot{\small\lfoot\ifnum\count0=\startpage\copyright\ \gtp\hfill\else
\agt, Volume \thevolumenumber\ (\thevolumeyear)\hfill\fi}
\def\@evenfoot{\@oddfoot}
\makeatother
\fi
\let\maketitlepage\makeagttitle
\let\makeshorttitle\maketitlepage
\let\maketitle\maketitlepage

\newwrite\gtoutfile
\long\gdef\makeheadfile{  
{\def\\{, }\def\s{ }
\immediate\openout\gtoutfile head.xxx
\immediate\write\gtoutfile{To: math@arxiv.org}
\immediate\write\gtoutfile{Subject: put OR rep NNNNN:ppppp}
\immediate\write\gtoutfile{--text follows this line--}
\immediate\write\gtoutfile{Proxy-for: \ifx\theasciiauthors\relax
\theauthors\else\theasciiauthors\fi\s<\ifx\theasciiemail\relax\theemail\else\theasciiemail\fi>}
\immediate\write\gtoutfile{\noexpand\\}
\immediate\write\gtoutfile{Authors: \ifx\theasciiauthors\relax
\theauthors\else\theasciiauthors\fi}
{\def\\{ }\immediate\write\gtoutfile{Title: \ifx\theasciititle\relax
\thetitle\else\theasciititle\fi}}
\immediate\write\gtoutfile{Subj-class: GT or SG, GR etc}
\immediate\write\gtoutfile{MSC-class: \theprimaryclass\ifx\thesecondaryclass\relax\else, \thesecondaryclass\fi}
\immediate\write\gtoutfile{Journal-ref: Algebr. Geom. Topol. \thevolumenumber\s
(\thevolumeyear) \startpage-\finishpage}
\immediate\write\gtoutfile{Comments: Published by Algebraic and
Geometric Topology at}
\immediate\write\gtoutfile{\s\s\s  http://www.maths.warwick.ac.uk/agt/AGTVol\thevolumenumber/agt-\thevolumenumber-\thepapernumber.abs.html}
\immediate\write\gtoutfile{\noexpand\\}
\immediate\write\gtoutfile{}
\ifx\theasciiabstract\relax
\immediate\write\gtoutfile{\theabstract}\else
\immediate\write\gtoutfile{\theasciiabstract}\fi
\immediate\write\gtoutfile{}
\immediate\write\gtoutfile{\noexpand\\}
\immediate\write\gtoutfile{}
\immediate\closeout\gtoutfile}}  

\def\maketitlepage{\makeagttitle\makeheadfile}
\let\makeshorttitle\maketitlepage
\let\maketitle\maketitlepage

\def\ifplaintex{\expandafter\ifx\csname documentclass\endcsname\relax}

\def\gtp{{\mathsurround=0pt\it $\cal G\mskip-2mu$eometry \&\ 
$\cal T\!\!$opology $\cal P\!$ublications}}  

\def\recd{{\small Received:\qua\receiveddate\ifx\reviseddate\relax
\else\qquad Revised:\qua\reviseddate\fi\par}} 

\def\lognumber#1{\def\thelognumber{#1}}
\def\volumenumber#1{\def\thevolumenumber{#1}}
\def\volumeyear#1{\def\thevolumeyear{#1}}
\def\papernumber#1{\def\thepapernumber{#1}}
\def\pagenumbers#1#2{\def\startpage{#1}\def\finishpage{#2}}
\def\published#1{\def\publishdate{#1}}

\def\received#1{\def\receiveddate{#1}}
\def\revised#1{\def\reviseddate{#1}}
\def\accepted#1{\def\accepteddate{#1}}
\def\asciititle#1{\def\theasciititle{#1}}
\def\covertitle#1{\def\thecovertitle{#1}}

\long\def\asciiabstract#1{\long\def\theasciiabstract{#1}}

\let\\\par\let\thelognumber\relax\let\thevolumenumber\relax
\let\thepapernumber\relax\let\thevolumeyear\relax\let\startpage\relax
\let\finishpage\relax\let\publishdate\relax\let\receiveddate\relax
\let\reviseddate\relax\let\accepteddate\relax\let\theasciititle\relax
\let\thecovertitle\relax\let\theasciiauthors\relax
\let\theasciiabstract\relax

\let\theasciiemail\relax

\ifplaintex
\font\logobig=cmssbx10 scaled 3836
\font\logomed=cmssbx10 scaled 2557
\else
\font\logobig=cmssbx10 scaled 4200
\font\logomed=cmssbx10 scaled 2800
\fi

\long\def\makeagttitle{   
\count0=\startpage
\agt\hfill      
\hbox to 45truept{\vbox to 0pt{\vglue -13truept{\logomed A\kern -.37em{\logobig 
T}\kern -.38em G}\vss}\hss}
\break
{\small Volume \thevolumenumber\ (\thevolumeyear)
\startpage--\finishpage\nl
Published: \publishdate}

\vglue .25truein

{\parskip=0pt\leftskip 0pt plus
1fil\def\\{\par\smallskip}{\Large\bf\thetitle}\par\medskip} \vglue
0.05truein

{\parskip=0pt\leftskip 0pt plus 1fil\def\\{\par}{\sc\theauthors}
\par\medskip}%
 
\vglue 0.03truein 

{\small\leftskip 25truept\rightskip 25truept{\bf Abstract}\stdspace\theabstract

{\bf AMS Classification}\stdspace\theprimaryclass
\ifx\thesecondaryclass\relax\else; \thesecondaryclass\fi\par
{\bf Keywords}\stdspace \thekeywords\par}\vglue 7truept

}   

\ifplaintex
\font\phead=cmsl9 scaled 950
\font\pnum=cmbx10 scaled 913
\font\pfoot=cmsl9 scaled 950
\headline{\vbox to 0pt{\vskip -4.5mm\line{\small\phead\ifnum
\count0=\startpage ISSN 1472-2739 (on-line) 1472-2747 (printed)
\hfill {\pnum\folio}\else\ifodd\count0\def\\{ }%
\ifx\theshorttitle\relax\thetitle\else\theshorttitle\fi\hfill{\pnum\folio}
\else\def\\{ and }{\pnum\folio}\hfill\ifx\theshortauthors\relax\theauthors
\else\theshortauthors\fi\fi\fi}\vss}}
\footline{\vbox to 0pt{\vglue 0mm\line{\small\pfoot\ifnum\count0=\startpage
\copyright\ \gtp\hfill\else
\agt, Volume \thevolumenumber\ (\thevolumeyear)\hfill\fi}\vss}}
\else
\font=cmsl9 scaled 1050
\font\lnum=cmbx10 
\font=cmsl9 scaled 1050
\makeatletter
\def\@oddhead{{\small\lhead\ifnum\count0=\startpage ISSN 1472-2739 
(on-line) 1472-2747 (printed)\hfill {\lnum\number\count0}\else\ifodd\count0
\def\\{ }\ifx\theshorttitle\relax \thetitle \else\theshorttitle\fi\hfill
{\lnum\number\count0}\else\def\\{ and }{\lnum\number\count0}
\hfill\ifx\theshortauthors\relax 
\theauthors\else\theshortauthors\fi\fi\fi}}\def\@evenhead{\@oddhead}
\def\@oddfoot{\small\lfoot\ifnum\count0=\startpage\copyright\ \gtp\hfill\else
\agt, Volume \thevolumenumber\ (\thevolumeyear)\hfill\fi}
\def\@evenfoot{\@oddfoot}
\makeatother
\fi
\let\maketitlepage\makeagttitle
\let\makeshorttitle\maketitlepage
\let\maketitle\maketitlepage

\newwrite\gtoutfile
\long\gdef\makeheadfile{  
{\def\\{, }\def\s{ }
\immediate\openout\gtoutfile head.xxx
\immediate\write\gtoutfile{To: math@arxiv.org}
\immediate\write\gtoutfile{Subject: put OR rep NNNNN:ppppp}
\immediate\write\gtoutfile{--text follows this line--}
\immediate\write\gtoutfile{Proxy-for: \ifx\theasciiauthors\relax
\theauthors\else\theasciiauthors\fi\s<\ifx\theasciiemail\relax\theemail\else\theasciiemail\fi>}
\immediate\write\gtoutfile{\noexpand\\}
\immediate\write\gtoutfile{Authors: \ifx\theasciiauthors\relax
\theauthors\else\theasciiauthors\fi}
{\def\\{ }\immediate\write\gtoutfile{Title: \ifx\theasciititle\relax
\thetitle\else\theasciititle\fi}}
\immediate\write\gtoutfile{Subj-class: GT or SG, GR etc}
\immediate\write\gtoutfile{MSC-class: \theprimaryclass\ifx\thesecondaryclass\relax\else, \thesecondaryclass\fi}
\immediate\write\gtoutfile{Journal-ref: Algebr. Geom. Topol. \thevolumenumber\s
(\thevolumeyear) \startpage-\finishpage}
\immediate\write\gtoutfile{Comments: Published by Algebraic and
Geometric Topology at}
\immediate\write\gtoutfile{\s\s\s  http://www.maths.warwick.ac.uk/agt/AGTVol\thevolumenumber/agt-\thevolumenumber-\thepapernumber.abs.html}
\immediate\write\gtoutfile{\noexpand\\}
\immediate\write\gtoutfile{}
\ifx\theasciiabstract\relax
\immediate\write\gtoutfile{\theabstract}\else
\immediate\write\gtoutfile{\theasciiabstract}\fi
\immediate\write\gtoutfile{}
\immediate\write\gtoutfile{\noexpand\\}
\immediate\write\gtoutfile{}
\immediate\closeout\gtoutfile}}  

\def\maketitlepage{\makeagttitle\makeheadfile}
\let\makeshorttitle\maketitlepage
\let\maketitle\maketitlepage

\def\ifplaintex{\expandafter\ifx\csname documentclass\endcsname\relax}

\def\gtp{{\mathsurround=0pt\it $\cal G\mskip-2mu$eometry \&\ 
$\cal T\!\!$opology $\cal P\!$ublications}}  

\def\recd{{\small Received:\qua\receiveddate\ifx\reviseddate\relax
\else\qquad Revised:\qua\reviseddate\fi\par}} 

\def\lognumber#1{\def\thelognumber{#1}}
\def\volumenumber#1{\def\thevolumenumber{#1}}
\def\volumeyear#1{\def\thevolumeyear{#1}}
\def\papernumber#1{\def\thepapernumber{#1}}
\def\pagenumbers#1#2{\def\startpage{#1}\def\finishpage{#2}}
\def\published#1{\def\publishdate{#1}}

\def\received#1{\def\receiveddate{#1}}
\def\revised#1{\def\reviseddate{#1}}
\def\accepted#1{\def\accepteddate{#1}}
\def\asciititle#1{\def\theasciititle{#1}}
\def\covertitle#1{\def\thecovertitle{#1}}

\long\def\asciiabstract#1{\long\def\theasciiabstract{#1}}

\let\\\par\let\thelognumber\relax\let\thevolumenumber\relax
\let\thepapernumber\relax\let\thevolumeyear\relax\let\startpage\relax
\let\finishpage\relax\let\publishdate\relax\let\receiveddate\relax
\let\reviseddate\relax\let\accepteddate\relax\let\theasciititle\relax
\let\thecovertitle\relax\let\theasciiauthors\relax
\let\theasciiabstract\relax

\let\theasciiemail\relax

\ifplaintex
\font\logobig=cmssbx10 scaled 3836
\font\logomed=cmssbx10 scaled 2557
\else
\font\logobig=cmssbx10 scaled 4200
\font\logomed=cmssbx10 scaled 2800
\fi

\long\def\makeagttitle{   
\count0=\startpage
\agt\hfill      
\hbox to 45truept{\vbox to 0pt{\vglue -13truept{\logomed A\kern -.37em{\logobig 
T}\kern -.38em G}\vss}\hss}
\break
{\small Volume \thevolumenumber\ (\thevolumeyear)
\startpage--\finishpage\nl
Published: \publishdate}

\vglue .25truein

{\parskip=0pt\leftskip 0pt plus
1fil\def\\{\par\smallskip}{\Large\bf\thetitle}\par\medskip} \vglue
0.05truein

{\parskip=0pt\leftskip 0pt plus 1fil\def\\{\par}{\sc\theauthors}
\par\medskip}%
 
\vglue 0.03truein 

{\small\leftskip 25truept\rightskip 25truept{\bf Abstract}\stdspace\theabstract

{\bf AMS Classification}\stdspace\theprimaryclass
\ifx\thesecondaryclass\relax\else; \thesecondaryclass\fi\par
{\bf Keywords}\stdspace \thekeywords\par}\vglue 7truept

}   

\ifplaintex
\font\phead=cmsl9 scaled 950
\font\pnum=cmbx10 scaled 913
\font\pfoot=cmsl9 scaled 950
\headline{\vbox to 0pt{\vskip -4.5mm\line{\small\phead\ifnum
\count0=\startpage ISSN 1472-2739 (on-line) 1472-2747 (printed)
\hfill {\pnum\folio}\else\ifodd\count0\def\\{ }%
\ifx\theshorttitle\relax\thetitle\else\theshorttitle\fi\hfill{\pnum\folio}
\else\def\\{ and }{\pnum\folio}\hfill\ifx\theshortauthors\relax\theauthors
\else\theshortauthors\fi\fi\fi}\vss}}
\footline{\vbox to 0pt{\vglue 0mm\line{\small\pfoot\ifnum\count0=\startpage
\copyright\ \gtp\hfill\else
\agt, Volume \thevolumenumber\ (\thevolumeyear)\hfill\fi}\vss}}
\else
\font=cmsl9 scaled 1050
\font\lnum=cmbx10 
\font=cmsl9 scaled 1050
\makeatletter
\def\@oddhead{{\small\lhead\ifnum\count0=\startpage ISSN 1472-2739 
(on-line) 1472-2747 (printed)\hfill {\lnum\number\count0}\else\ifodd\count0
\def\\{ }\ifx\theshorttitle\relax \thetitle \else\theshorttitle\fi\hfill
{\lnum\number\count0}\else\def\\{ and }{\lnum\number\count0}
\hfill\ifx\theshortauthors\relax 
\theauthors\else\theshortauthors\fi\fi\fi}}\def\@evenhead{\@oddhead}
\def\@oddfoot{\small\lfoot\ifnum\count0=\startpage\copyright\ \gtp\hfill\else
\agt, Volume \thevolumenumber\ (\thevolumeyear)\hfill\fi}
\def\@evenfoot{\@oddfoot}
\makeatother
\fi
\let\maketitlepage\makeagttitle
\let\makeshorttitle\maketitlepage
\let\maketitle\maketitlepage

\newwrite\gtoutfile
\long\gdef\makeheadfile{  
{\def\\{, }\def\s{ }
\immediate\openout\gtoutfile head.xxx
\immediate\write\gtoutfile{To: math@arxiv.org}
\immediate\write\gtoutfile{Subject: put OR rep NNNNN:ppppp}
\immediate\write\gtoutfile{--text follows this line--}
\immediate\write\gtoutfile{Proxy-for: \ifx\theasciiauthors\relax
\theauthors\else\theasciiauthors\fi\s<\ifx\theasciiemail\relax\theemail\else\theasciiemail\fi>}
\immediate\write\gtoutfile{\noexpand\\}
\immediate\write\gtoutfile{Authors: \ifx\theasciiauthors\relax
\theauthors\else\theasciiauthors\fi}
{\def\\{ }\immediate\write\gtoutfile{Title: \ifx\theasciititle\relax
\thetitle\else\theasciititle\fi}}
\immediate\write\gtoutfile{Subj-class: GT or SG, GR etc}
\immediate\write\gtoutfile{MSC-class: \theprimaryclass\ifx\thesecondaryclass\relax\else, \thesecondaryclass\fi}
\immediate\write\gtoutfile{Journal-ref: Algebr. Geom. Topol. \thevolumenumber\s
(\thevolumeyear) \startpage-\finishpage}
\immediate\write\gtoutfile{Comments: Published by Algebraic and
Geometric Topology at}
\immediate\write\gtoutfile{\s\s\s  http://www.maths.warwick.ac.uk/agt/AGTVol\thevolumenumber/agt-\thevolumenumber-\thepapernumber.abs.html}
\immediate\write\gtoutfile{\noexpand\\}
\immediate\write\gtoutfile{}
\ifx\theasciiabstract\relax
\immediate\write\gtoutfile{\theabstract}\else
\immediate\write\gtoutfile{\theasciiabstract}\fi
\immediate\write\gtoutfile{}
\immediate\write\gtoutfile{\noexpand\\}
\immediate\write\gtoutfile{}
\immediate\closeout\gtoutfile}}  

\def\maketitlepage{\makeagttitle\makeheadfile}
\let\makeshorttitle\maketitlepage
\let\maketitle\maketitlepage

\def\ifplaintex{\expandafter\ifx\csname documentclass\endcsname\relax}

\def\gtp{{\mathsurround=0pt\it $\cal G\mskip-2mu$eometry \&\ 
$\cal T\!\!$opology $\cal P\!$ublications}}  

\def\recd{{\small Received:\qua\receiveddate\ifx\reviseddate\relax
\else\qquad Revised:\qua\reviseddate\fi\par}} 

\def\lognumber#1{\def\thelognumber{#1}}
\def\volumenumber#1{\def\thevolumenumber{#1}}
\def\volumeyear#1{\def\thevolumeyear{#1}}
\def\papernumber#1{\def\thepapernumber{#1}}
\def\pagenumbers#1#2{\def\startpage{#1}\def\finishpage{#2}}
\def\published#1{\def\publishdate{#1}}

\def\received#1{\def\receiveddate{#1}}
\def\revised#1{\def\reviseddate{#1}}
\def\accepted#1{\def\accepteddate{#1}}
\def\asciititle#1{\def\theasciititle{#1}}
\def\covertitle#1{\def\thecovertitle{#1}}

\long\def\asciiabstract#1{\long\def\theasciiabstract{#1}}

\let\\\par\let\thelognumber\relax\let\thevolumenumber\relax
\let\thepapernumber\relax\let\thevolumeyear\relax\let\startpage\relax
\let\finishpage\relax\let\publishdate\relax\let\receiveddate\relax
\let\reviseddate\relax\let\accepteddate\relax\let\theasciititle\relax
\let\thecovertitle\relax\let\theasciiauthors\relax
\let\theasciiabstract\relax

\let\theasciiemail\relax

\ifplaintex
\font\logobig=cmssbx10 scaled 3836
\font\logomed=cmssbx10 scaled 2557
\else
\font\logobig=cmssbx10 scaled 4200
\font\logomed=cmssbx10 scaled 2800
\fi

\long\def\makeagttitle{   
\count0=\startpage
\agt\hfill      
\hbox to 45truept{\vbox to 0pt{\vglue -13truept{\logomed A\kern -.37em{\logobig 
T}\kern -.38em G}\vss}\hss}
\break
{\small Volume \thevolumenumber\ (\thevolumeyear)
\startpage--\finishpage\nl
Published: \publishdate}

\vglue .25truein

{\parskip=0pt\leftskip 0pt plus
1fil\def\\{\par\smallskip}{\Large\bf\thetitle}\par\medskip} \vglue
0.05truein

{\parskip=0pt\leftskip 0pt plus 1fil\def\\{\par}{\sc\theauthors}
\par\medskip}%
 
\vglue 0.03truein 

{\small\leftskip 25truept\rightskip 25truept{\bf Abstract}\stdspace\theabstract

{\bf AMS Classification}\stdspace\theprimaryclass
\ifx\thesecondaryclass\relax\else; \thesecondaryclass\fi\par
{\bf Keywords}\stdspace \thekeywords\par}\vglue 7truept

}   

\ifplaintex
\font\phead=cmsl9 scaled 950
\font\pnum=cmbx10 scaled 913
\font\pfoot=cmsl9 scaled 950
\headline{\vbox to 0pt{\vskip -4.5mm\line{\small\phead\ifnum
\count0=\startpage ISSN 1472-2739 (on-line) 1472-2747 (printed)
\hfill {\pnum\folio}\else\ifodd\count0\def\\{ }%
\ifx\theshorttitle\relax\thetitle\else\theshorttitle\fi\hfill{\pnum\folio}
\else\def\\{ and }{\pnum\folio}\hfill\ifx\theshortauthors\relax\theauthors
\else\theshortauthors\fi\fi\fi}\vss}}
\footline{\vbox to 0pt{\vglue 0mm\line{\small\pfoot\ifnum\count0=\startpage
\copyright\ \gtp\hfill\else
\agt, Volume \thevolumenumber\ (\thevolumeyear)\hfill\fi}\vss}}
\else
\font=cmsl9 scaled 1050
\font\lnum=cmbx10 
\font=cmsl9 scaled 1050
\makeatletter
\def\@oddhead{{\small\lhead\ifnum\count0=\startpage ISSN 1472-2739 
(on-line) 1472-2747 (printed)\hfill {\lnum\number\count0}\else\ifodd\count0
\def\\{ }\ifx\theshorttitle\relax \thetitle \else\theshorttitle\fi\hfill
{\lnum\number\count0}\else\def\\{ and }{\lnum\number\count0}
\hfill\ifx\theshortauthors\relax 
\theauthors\else\theshortauthors\fi\fi\fi}}\def\@evenhead{\@oddhead}
\def\@oddfoot{\small\lfoot\ifnum\count0=\startpage\copyright\ \gtp\hfill\else
\agt, Volume \thevolumenumber\ (\thevolumeyear)\hfill\fi}
\def\@evenfoot{\@oddfoot}
\makeatother
\fi
\let\maketitlepage\makeagttitle
\let\makeshorttitle\maketitlepage
\let\maketitle\maketitlepage

\newwrite\gtoutfile
\long\gdef\makeheadfile{  
{\def\\{, }\def\s{ }
\immediate\openout\gtoutfile head.xxx
\immediate\write\gtoutfile{To: math@arxiv.org}
\immediate\write\gtoutfile{Subject: put OR rep NNNNN:ppppp}
\immediate\write\gtoutfile{--text follows this line--}
\immediate\write\gtoutfile{Proxy-for: \ifx\theasciiauthors\relax
\theauthors\else\theasciiauthors\fi\s<\ifx\theasciiemail\relax\theemail\else\theasciiemail\fi>}
\immediate\write\gtoutfile{\noexpand\\}
\immediate\write\gtoutfile{Authors: \ifx\theasciiauthors\relax
\theauthors\else\theasciiauthors\fi}
{\def\\{ }\immediate\write\gtoutfile{Title: \ifx\theasciititle\relax
\thetitle\else\theasciititle\fi}}
\immediate\write\gtoutfile{Subj-class: GT or SG, GR etc}
\immediate\write\gtoutfile{MSC-class: \theprimaryclass\ifx\thesecondaryclass\relax\else, \thesecondaryclass\fi}
\immediate\write\gtoutfile{Journal-ref: Algebr. Geom. Topol. \thevolumenumber\s
(\thevolumeyear) \startpage-\finishpage}
\immediate\write\gtoutfile{Comments: Published by Algebraic and
Geometric Topology at}
\immediate\write\gtoutfile{\s\s\s  http://www.maths.warwick.ac.uk/agt/AGTVol\thevolumenumber/agt-\thevolumenumber-\thepapernumber.abs.html}
\immediate\write\gtoutfile{\noexpand\\}
\immediate\write\gtoutfile{}
\ifx\theasciiabstract\relax
\immediate\write\gtoutfile{\theabstract}\else
\immediate\write\gtoutfile{\theasciiabstract}\fi
\immediate\write\gtoutfile{}
\immediate\write\gtoutfile{\noexpand\\}
\immediate\write\gtoutfile{}
\immediate\closeout\gtoutfile}}  

\def\maketitlepage{\makeagttitle\makeheadfile}
\let\makeshorttitle\maketitlepage
\let\maketitle\maketitlepage

\lognumber{19}
\volumenumber{1}
\volumeyear{2001}
\papernumber{19}
\pagenumbers{381}{409}
\received{7 December 2000}
\revised{24 May 2000}
\accepted{18 June 2001}
\published{19 June 2001}

\input xy
\xyoption{all}

\newtheorem{thrm}{Theorem}[section]
\newtheorem{lem}[thrm]{Lemma}
\newtheorem{cor}[thrm]{Corollary}
\newtheorem{prop}[thrm]{Proposition}

\theoremstyle{definition}

\newtheorem{rem}[thrm]{Remark}
\newtheorem{defn}[thrm]{Definition}
\newtheorem{exmpl}[thrm]{Example}

\def\maprt#1{\, \smash{\mathop{\longrightarrow}\limits^{#1}}\, }

\def\bd{  \begin{diagram}    }
\def\ed{  \end{diagram}      }

\def\thm #1  {\medskip\noindent{\bf #1}\qua}
\def\pf #1 {\proof[#1]}
\def\term #1{{\sl #1}}

\let\halmos\endproof

\def\Bbb{\bf}

\def\ZZ{{\Bbb Z}}

\def\HH{{\Bbb H}}

\def\RR{{\Bbb R}}

\def\CC{{\Bbb C}}

\def\FF{{\Bbb F}}

\def\s{\Sigma}
\def\smsh{\wedge}
\def\om{\Omega}

\def\cross{\times}
\def\wdg{\vee}

\def\inclds{\hookrightarrow}
\def\dim{{\rm dim}}

\def\of{\circ}
\def\twdl{\widetilde}
\def\<{\langle}
\def\>{\rangle}

\def\sseq{\subseteq}

\def\bprp{\begin{prop}}
\def\eprp{\end{prop}}

\def\bthm{\begin{thrm}}
\def\ethm{\end{thrm}}

\def\blem{\begin{lem}}
\def\elem{\end{lem}}

\def\bcor{\begin{cor}}
\def\ecor{\end{cor}}

\def\brmk{\begin{rem}}
\def\ermk{\end{rem}}

\def\bdfn{\begin{defn}}
\def\edfn{\end{defn}}

\def\bexm{\begin{exmpl}}
\def\eexm{\end{exmpl}}

\def\Im{ {\mathrm I}{\mathrm m} }
\def\im{\Im}

\def\Hom{{\mathrm H}{\mathrm o}{\mathrm m} }

\def\Im{\textup{Im}}

\def\Hom{\textup{Hom}}

\def\blorp#1{}

\def\lcm{\textup{lcm}}

\def\|{\, \bigm|\, }

\def\F{{\cal F}}
\def\M{{\cal M}}
\def\S{{\cal S}}
\def\E{{\cal E}}

\def\cp{\CC {\mathrm P} }
\def\rp{\RR {\mathrm P} }
\def\hp{\HH {\mathrm P} }
\def\fp{\FF {\mathrm P} }

\def\kl{{\mathrm{kl}}}
\def\cl{{\mathrm{cl}}}
\def\cat{{\mathrm{cat}}}

\def\Ext{{\mathrm{Ext}}}
\def\Hom{{\mathrm{Hom}}}

\def\simeq{=}

\begin{document}

\title{Homotopy classes that are trivial mod $\F$}
\asciititle{Homotopy classes that are trivial mod F}
\covertitle{Homotopy classes that are trivial mod $\noexpand\cal F$}
\authors{Martin Arkowitz and Jeffrey Strom}
\address{Dartmouth College, Hanover NH 03755, USA}
\email{martin.arkowitz@dartmouth.edu, jeffrey.strom@dartmouth.edu}

\begin{abstract}
If $\F$  is a collection of topological spaces,  then a 
homotopy class  $\alpha$  in  $[X,Y]$  is called   $\F$-trivial if   
$$
\alpha _* = 0: [A,X] \maprt{} [A,Y]
$$
for all  $A\in\F$.  In this paper we study the collection  
$Z_{\F}(X,Y)$  of all  $\F$-trivial  homotopy classes in  $[X,Y]$  
when  $\F = \S$,  the collection of spheres,  $\F = \M$,  the 
collection of Moore spaces, and  $F = \Sigma$,  the collection 
of suspensions.  Clearly  
$$
Z_{\Sigma}(X,Y) \subseteq Z_{\M}(X,Y) \subseteq Z_{\S}(X,Y),
$$ 
and we find examples of {\it finite complexes}  $X$  and  $Y$  for 
which these inclusions are strict.  We are also interested in  
$Z_{\F}(X) = Z_{\F}(X,X)$,  which under composition has the 
structure of a semigroup with zero.  We show that if  $X$  is a  
finite dimensional complex and  $\F =\S$,  $\M$  or  $\Sigma$,  then 
the semigroup  $Z_{\F}(X)$ is nilpotent.  More precisely,  the 
nilpotency of  $Z_{\F}(X)$  is bounded above by the  $\F$-killing  
length of  $X$,  a new numerical invariant which equals the 
number of steps it takes to make  $X$  contractible by 
successively attaching cones on wedges of spaces in  $\F$,  and 
this in turn is bounded above by the  $\F$-cone length of  X.  We 
then calculate or estimate the nilpotency of  $Z_{\F}(X)$  when  
$\F = \S$, $\M$  or $\Sigma$  for the following classes of spaces: 
(1) projective spaces  (2) certain Lie groups such as  $SU(n)$ 
and  $Sp(n)$.  The paper concludes with several open 
problems.
\end{abstract}

\asciiabstract{If F is a collection of topological spaces, then a
homotopy class \alpha in [X,Y] is called F-trivial if \alpha _* = 0:
[A,X] --> [A,Y] for all A in F.  In this paper we study the collection
Z_F(X,Y) of all F-trivial homotopy classes in [X,Y] when F = S, the
collection of spheres, F = M, the collection of Moore spaces, and F =
\Sigma, the collection of suspensions.  Clearly Z_\Sigma(X,Y) \subset
Z_M(X,Y) \subset Z_S(X,Y), and we find examples of finite complexes X
and Y for which these inclusions are strict.  We are also interested
in Z_F(X) = Z_F(X,X), which under composition has the structure of a
semigroup with zero.  We show that if X is a finite dimensional
complex and F = S, M or \Sigma, then the semigroup Z_F(X) is
nilpotent.  More precisely, the nilpotency of Z_F(X) is bounded above
by the F-killing length of X, a new numerical invariant which equals
the number of steps it takes to make X contractible by successively
attaching cones on wedges of spaces in F, and this in turn is bounded
above by the F-cone length of X.  We then calculate or estimate the
nilpotency of Z_F(X) when F = S, M or \Sigma for the following classes
of spaces: (1) projective spaces (2) certain Lie groups such as
SU(n) and Sp(n).  The paper concludes with several open problems.}

\primaryclass{55Q05}
\secondaryclass{55P65, 55P45, 55M30}
\keywords{Cone length, trivial homotopy}                    

\maketitle

\section{Introduction}

A map $f:X\maprt{}Y$ is said to be detected by
a collection $\F$ of topological spaces if there is a 
space $A\in\F$ such that the induced map 
$f_*:[A,X]\maprt{} [A,Y]$ of homotopy sets is nontrivial.  
It is a standard technique in homotopy theory to use certain simple 
collections $\F$
to detect essential homotopy classes.  
In studying the entire homotopy set $[X,Y]$ using this approach, one is led 
naturally
to consider the set of homotopy classes which are {\it not} detected
by $\F$, called the $\F$-trivial homotopy classes.
For example, if $\S$ is the collection of spheres, then
$f:X\maprt{}Y$ is detected by $\S$ precisely when some
induced homomorphism of homotopy groups $\pi_k(f):\pi_k(X)\maprt{}\pi_k(Y)$
is nonzero.  
The $\S$-trivial homotopy classes are those that induce zero on all homotopy groups.  
It is also important to determine induced maps on homotopy sets.
For this, one needs to understand composition of $\F$-trivial homotopy
classes.
With this in mind, we study two basic questions in this paper 
for a fixed collection $\F$: (1)  What is the 
set of all $\F$-trivial homotopy classes in $[X,Y]$? and 
(2) In the special case $X=Y$, how do $\F$-trivial homotopy classes behave 
under composition?   We are particularly interested in the
collections $\S$ of spheres, $\M$ of Moore spaces and $\s$ of suspensions.

Some of these ideas have appeared earlier.  The paper \cite{AMS}
considers the special case $\F=\S$.  Furthermore, Christensen has
studied similar questions in the stable category \cite{Christensen}.

We next briefly summarize the contents of this paper.  We write $Z_\F(X,Y)$
for the $\F$-trivial homotopy classes in $[X,Y]$ and set $Z_\F(X)=Z_\F(X,X)$.
After some generalities on $Z_\F(X,Y)$, we observe in Section 2 that $Z_\F(X)$
is a semigroup under composition.  Its nilpotency, denoted $t_\F(X)$, is a new
numerical invariant of homotopy type.  For the collection of suspensions, 
we prove that $t_\s(X)\leq \left\lceil \log_2(\cat(X))\right\rceil $.  In 
Section
3 we relate $t_\F(X)$ to other numerical invariants for arbitrary collections 
$\F$.
The $\F$-killing length of $X$, denoted $\kl_\F(X)$ (resp., the $\F$-cone 
length
of $X$, denoted $\cl_\F(X)$), is the least number of steps needed to go from 
$X$ to a contractible space (resp., from a contractible space to $X$) by successively
attaching cones on wedges of spaces in $\F$.  We prove that $t_\F(X)\leq \kl_\F(X)$,
and, if $\F$ is closed under suspension, that $\kl_\F(X)\leq \cl_\F(X)$.  We also
show that $\kl_\F(X)$ behaves subadditively with respect to cofibrations.
It is clear that for any $X$ and $Y$, $Z_\s(X,Y)\sseq Z_\M(X,Y)\sseq Z_\S(X,Y)$,
and we ask in Section 4 if these containments can be strict.  
It is easy to find 
infinite complexes with strict containment.  However, in Section 4 we solve
the more difficult problem of finding finite complexes with this property.  From 
this, we deduce that
containment can be strict for finite complexes when 
$X=Y$.  The next two sections are devoted to determining $Z_\F(X)$
and $t_\F(X)$ for certain classes of spaces.  In Section 5 we 
calculate $Z_\F(X)$ and $t_\F(X)$ for $\F=\S, \M$ and $\s$ when $X$
is any real or complex projective space, or is the  quaternionic projective
space $\hp^n$ with $n\leq 4$.  In Section 6 we consider $t_\s(Y)$ for 
certain Lie groups $Y$.  We show that $2\leq t_\s(Y)$ when $Y=SU(n)$ or
$Sp(n)$ by proving that the groups $[Y,Y]$ are not abelian.  In addition,
we compute $t_\s(SO(n))$ for $n=3$ and $4$.  The paper concludes in Section
7 with a list of open problems.

For the remainder of this section, we give our notation and terminology.
All topological spaces are based and connected, and have the based homotopy type
of CW complexes.  All maps and homotopies preserve base points.   
We do not usually distinguish notationally between a map and its homotopy class.   
We let $*$ denote  the base point of a space or a space consisting of a single point.
In addition to standard notation, we use
$\equiv$ for same homotopy type,  $0\in [X,Y]$ for the constant homotopy class 
and $\mathrm{id}\in[X,X]$ for
the identity homotopy class.

For an abelian group $G$ and an integer $n\geq 2$, we let $M(G,n)$ denote the Moore
space of type $(G,n)$, that is, the space with a single non-vanishing reduced homology
group $G$ in dimension $n$. If $G$ is finitely generated, we also define $M(G,1)$
as a wedge of circles $S^1$ and spaces obtained by 
attaching a 2-cell to $S^1$ by a map of degree $m$.  The $n^{th}$ homotopy group of 
$X$ with coefficients
in $G$ is $\pi_n(X;G)= [M(G,n), X]$.  A map $f:X\maprt{}Y$ induces a homomorphism
$\pi_n(f;G):\pi_n(X;G)\maprt{} \pi_n(Y;G)$, and $\pi_*(f;G)$ denotes the set of 
all such homomorphisms.  If $G=\ZZ$, we write $\pi_n(X)$ and $\pi_n(f)$ for the 
$n^{th}$ homotopy group and induced map, respectively.

We use unreduced Lusternik-Schnirelmann category of a space $X$, denoted $\cat(X)$.
Thus $\cat(X)\leq 2$ if and only if $X$ is a co-H-space.
By an H-space, we mean a space with a homotopy-associative multiplication and 
homotopy inverse, i.e., a group-like space.

For a positive integer $n$, the cyclic group of order $n$ is denoted $\ZZ/n$.
If $X$ is a space or an abelian group, we use the notation $X_{(p)}$ for the localization
of $X$ at the prime $p$ \cite{HMR}.  We let $\lambda: X\maprt{} X_{(p)}$ denote the natural map
from $X$ to its localization.

A semigroup is a set $S$ with an associative binary operation, denoted by juxtaposition.
We call $S$ a pointed semigroup if there is an element $0\in S$ such that 
$x0 =0x =0$ for each $x\in S$.  A pointed semigroup is nilpotent if there
is an integer $n$ such that the product
$x_1\cdots x_n$ is $0$ whenever $x_1, \ldots , x_n\in S$.  
The least such integer $n$ is the nilpotency of $S$.
If $S$ is not nilpotent, 
then we say its nilpotency is
$\infty$.  Finally, if $x$ is a real number, then $\lceil x
\rceil$ denotes the least integer $n\geq x$.

\section{$\F$-trivial homotopy classes}

Let $\F$ be any collection of spaces.  

\bdfn
A homotopy class $f:X\maprt{}Y$ is 
\term{$\F$-trivial} if the induced map $f_*:[A,X]\maprt{} [A,Y]$
is trivial for each $A\in \F$.  
We denote by $Z_\F (X,Y)$ the subset of $[X,Y]$ consisting of all
$\F$-trivial homotopy classes.  We denote
$Z_\F(X,X)$ by $Z_\F(X)$. 
\edfn

We study $Z_\F(X,Y)$ and $Z_\F(X)$ for 
certain collections $\F$.  The following are some
interesting examples.

\thm{Examples}
\begin{enumerate}
\item[(a)] $\S=\{ S^n \| n\geq 1\}$, the collection of spheres.  
In this case $f\in Z_\S(X,Y)$ if and only
if $\pi_*(f) =0$.
\item[(b)]  $\M=\{ M(\ZZ/m,n) \|  m\geq 0, n\geq 1\}$, 
the collection of Moore spaces $M(\ZZ/m,n)$.  Here $f\in Z_\M(X,Y)$
if and only if $\pi_*(f;G)=0$  for  any 
finitely generated abelian group $G$.
\item[(c)]  $\Sigma=\{ \s A \}$, the collection of all suspensions.  
In this case $f\in Z_\Sigma(X,Y)$ if and
only if $f_*=0:[\s A, X]\maprt{}[\s A, Y]$ for every space $A$.
\item[(d)]  ${\cal P}$ is the collection of all finite dimensional 
complexes.\ Then
$f\in Z_{\cal P}(X,Y)$ if and only if $f:X\maprt{} Y$ is a phantom map \cite{McG}.
\end{enumerate}

In this paper our main interest is in the collections $\S$, $\M$ and $\s$.  
In Section 7 we will mention a few other collections.
However, we begin with some simple, general facts about arbitrary collections.

\blem\label{lem:basic}\hfill
\begin{enumerate}
\item[\rm(a)]
If $\F\sseq \F'$, then $Z_{\F'}(X,Y) \sseq Z_{\F}(X,Y)$ for any 
$X$ and $Y$.  
\item[\rm(b)]  If $X_\alpha \in \F$ for each $\alpha$ in some index set, 
then $Z_\F(\bigvee X_\alpha ,Y)=  0 $ for every $Y$. 
\item[\rm(c)]  For any $X$,  $Z_\F(X)$ is a pointed semigroup under 
the binary operation of composition of homotopy classes, and with 
zero the constant homotopy class.
\end{enumerate}
\elem

Any map $f:X\maprt{}Y$ gives rise to functions $\twdl f:X\cross Y\maprt{} X\cross Y$
and $\overline f: X\wdg Y \maprt{} X\wdg Y$ defined by the 
diagrams
$$
\xymatrix{
X\cross Y \ar[r]^{\twdl f}\ar[d]_{p_1} & X\cross Y\\
X\ar[r]^f & Y\ar[u]_{i_2}\\ }
\qquad
\mathrm{and}
\qquad
\xymatrix{
X\wdg Y \ar[r]^{\overline f}\ar[d]_{q_1} & X\wdg Y\\
X\ar[r]^f & Y.\ar[u]_{j_2}\\ }
$$
Clearly, if $f\in Z_\F(X,Y)$, then $\twdl f\in Z_\F(X\cross Y)$ and $\overline f\in Z_\F(X\wdg Y)$.

The following lemma, whose proof is obvious, will be used frequently.

\blem\label{lem:self}
The functions 
$$
\theta:Z_\F(X,Y) \maprt{} Z_\F(X\cross Y)
\quad \mathrm{and}\quad \phi:Z_\F(X,Y) \maprt{} Z_\F(X\wdg Y),
$$ 
defined by $\theta(f) = \twdl f$ and $\phi(f) = \overline f$, are 
injective.  Thus, $Z_\F(X,Y) \neq 0$ implies that
$Z_\F(X\cross Y)\neq 0$ and $Z_\F(X\wdg Y)\neq 0$
\elem

We conclude this section with some basic
results about $Z_\Sigma(X,Y)$ and $Z_\Sigma(X)$.
Recall that a map $f:X\maprt{}Y$ has \term{essential category weight}
at least $n$, written $E(f)\geq n$, if for every space $A$ with $\cat(A)\leq n$,
we have $f_*=0 :[A,X]\maprt{} [A,Y]$ \cite{Strom,Rudyak}.

\blem\label{lem:E}  For any two spaces $X$ and $Y$, 
$$
\begin{array}[]{ll}
Z_\s(X,Y) &  =  \{ f\, |\, f\in [X,Y],\, E(f)\geq 2     \} \\
          &  =  \{ f\, |\, f\in [X,Y],\,  \om f\simeq 0  \} .
\end{array}
$$
\elem

\pf Proof
Let $f\in Z_\s(X,Y)$.  If $\cat(A)\leq 2$ then 
the canonical map $\nu\co\s\om A\maprt{}A$ has a section $s$.  Thus the diagram
$$
\xymatrix{
[A,X] \ar[d]_{\nu^*}
      \ar[r]^{f_{*}}     & 
[A,Y]\ar[d]^{\nu^*}      \\
[\s\om A, X]
       \ar[r]^{f_{**}}      & 
[\s\om A, Y]             \\          }
$$
commutes, and so $f_* = s^* f_{**}\nu^* = 0$.  Since the reverse implication is trivial,
this establishes the first equality.

Now assume that $f_{*}=0: [\s B, X]\maprt{}[\s B, Y]$
for every space $B$.  Taking $B=\om X$, we find that $f\of\nu \simeq 0: \s\om X\maprt{}Y$.
Since this map is adjoint to $\om f$, we conclude that $\om f\simeq 0$.  Conversely, if 
$\om f\simeq 0$, then $\om f_* = 0:[ B,\om X]\maprt{} [B,\om Y]$, which means that
$f_*=0 : [\s B, X]\maprt{} [\s B, Y]$.  This completes the proof.
\halmos

\thm{Remarks}
\begin{enumerate}
\item[(a)]  Since $\cat(A)\leq 2$ if and only if $A$ is a
co-H-space, a map $f:X\maprt{}Y$ has $E(f)\geq 2$ if and only if 
$f_*=0 : [A,X]\maprt{} [A,Y]$
for every co-H-space $A$. 
\item[(b)]
By Lemma \ref{lem:E}, we can regard the set $Z_\s(X,Y)$ as the kernel of the looping
function $\om : [X,Y]\maprt{} [\om X,\om Y]$.  We see from (a) that $\ker\om =0$ if 
$X$ is a co-H-space.  The function $\om$ has been extensively studied
in special cases, e.g., when $Y$ is an Eilenberg-MacLane space, then $\om$ is just the 
cohomology suspension 
\cite[Chap.$\,$VII]{Wh}.
\end{enumerate}

\bprp\label{prop:Et}
Let $X$ be a space of finite category, and let  $n\!\geq \!\log_2(\cat(X))$.
If $f_1,\ldots, f_n\in Z_\s(X)$, then $f_1\of\cdots \of f_n =0$.
Thus the nilpotency of the semigroup $Z_\s(X)$ is at most 
$\left\lceil \log_2(\cat(X))\right\rceil $, the least integer greater 
than or equal to $\log_2(\cat(X))$.
\eprp

\pf Proof
Since $f_i\in Z_\s(X)$, Lemma \ref{lem:E} shows that $E(f_i)\geq 2$.  
By the product formula for essential category weight  \cite[Thm.$\,$9]{Strom}, 
$E(f_1\of \cdots\of f_n)\geq E(f_1)\cdots E(f_n) \geq 2^n \geq \cat(X)$. From 
the definition of essential category weight, 
$f_1\of \cdots\of f_n=0$.\halmos

\thm{Remark} 
We shall see later  that 
the semigroup $Z_\S(X)$ is nilpotent if $X$ is a finite
dimensional complex.  It follows that this is true for $Z_\M(X)$
and $Z_\s(X)$  (Remark (b) following Theorem
\ref{thrm:nilbound}).

\bdfn
For any collection $\F$ of spaces and any space $X$, we define
$t_\F(X)$, the \term{nilpotency of $X$ mod $\F$} as follows: If $X$ is contractible, 
set $t_\F(X)=0$;  Otherwise, $t_\F(X)$ is the nilpotency of the semigroup $Z_\F(X)$.
\edfn

Thus $t_\F(X) = 1$ if and only if $X$ is not contractible and $Z_\F(X) =0$.

The set $Z_\S(X)$ and the integer $t_\S(X)$ were considered in \cite{AMS}, where they
were written $Z_\infty(X)$ and $t_\infty(X)$.   Since $\S\sseq \M\sseq \s$, we have 
$$
0\leq t_\s (X) \leq t_\M (X) \leq t_\S(X) \leq \infty
$$
for any space $X$.

Since $\cat( \s A_1\cross \cdots \cross \s A_r)\leq r+1$ 
\cite[Prop.$\,$2.3]{James}, 
we have the following result.

\bcor\label{cor:spheret}
For any $r$ spaces $A_1,\ldots, A_r$, 
$$
t_\Sigma (\s A_1\cross\cdots\cross \s A_r)\leq \big\lceil \log_2( r+1)\big\rceil .
$$
\ecor

This paper is devoted to a study of the
sets $Z_\F(X,Y)$, with emphasis on the  
nilpotency of spaces mod $\F$ for $\F =\S,\M$ and $\s$.

\section{$\F$-killing length and $\F$-cone length}

Proposition \ref{prop:Et} shows that
$\lceil \log_2( \cat(X) )\rceil$ is an upper bound for 
$t_\s(X)$.  In this section, we obtain upper bounds on 
$t_\F(X)$ for arbitrary collections $\F$.
We begin with the main definitions of this section.

\bdfn\label{defn:killcone}
Let $\F$ be a collection of spaces and $X$ a space.  Suppose there is a
sequence of cofibrations
$$
L_i\maprt{} X_i\maprt{} X_{i+1}
$$
for $0\leq i < m$
such that each $L_i$ is a wedge of spaces which belong to $\F$.
If $X_0 \equiv X$ and $X_m\equiv *$, then this 
is called an \term{$\F$-killing length decomposition} of $X$ with length $m$. 
If $X_0\equiv *$ and $X_m\equiv X$,
then this is an \term{$\F$-cone length decomposition} with length $m$.
Define the \term{$\F$-killing length} and the \term{$\F$-cone length} 
of $X$, denoted by $\kl_\F(X)$ and $\cl_\F(X)$, respectively, as follows.
If $X\equiv *$, then $\kl_\F(X)=0$; otherwise, $\kl_\F(X)$ is the smallest integer
$m$ such that there exists an $\F$-killing length decomposition of $X$ with length
$m$.  The $\F$-cone length of $X$ is defined analogously.
\edfn

The main result of this section is that $\kl_\F(X)$ is an upper bound for 
$t_\F(X)$. We need a lemma.

\blem\label{lem:coflem}
If $X\maprt{f}Y\maprt{g}Z$ is a sequence of spaces and maps,
then there is a cofiber sequence of mapping cones
$C_f \maprt{} C_{g f}\maprt{} C_g$, where the maps are induced by $f$ and $g$.
\elem

The proof is elementary, and hence omitted.

\bthm\label{thrm:nilbound}
If $\F$ is any collection of spaces and $X$ is any space, then 
$$
t_\F(X)\leq \kl_\F(X).
$$
If $\F$ is closed under suspensions, then $\kl_\F(X)\leq \cl_\F(X)$.
\ethm

\pf Proof
Assume that $\kl_\F(X)=m >0$ with $\F$-killing length decomposition
$$
L_i\maprt{f_i} X_i\maprt{p_i} X_{i+1}
$$
for $0\leq i<m$.  Let
$g_0, \ldots , g_{m-1}\in Z_\F(X)$ and 
consider the following diagram, with dashed arrows to be inductively defined below:\eject

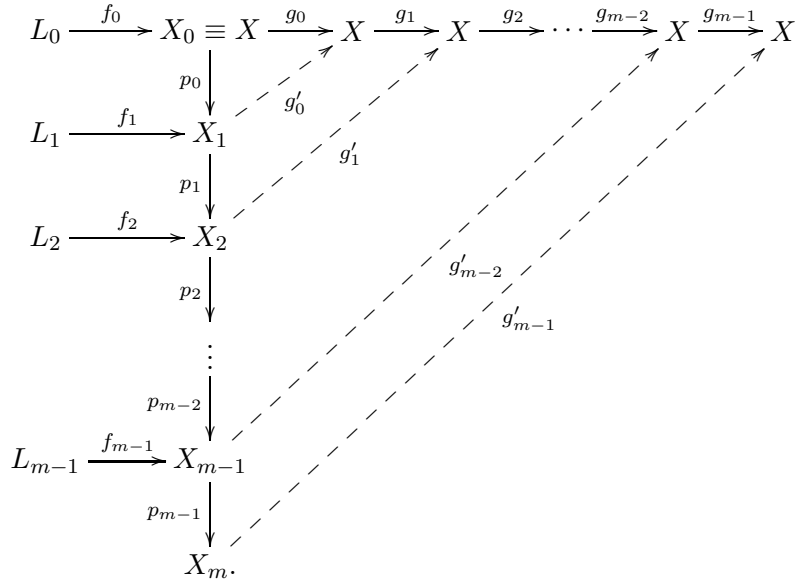
$$
\xymatrix{
L_{0}\ar[r]^(.4){f_0}& 
  X_{0}\equiv X \ar[d]_{p_{0}}\ar[r]^(.6){g_0} & 
  X \ar[r]^{g_1} &
  X \ar[r]^{g_2} &
  \cdots \ar[r]^{g_{m-2}} &
  X \ar[r]^{g_{m-1}} &   
  X  \\
L_{1}\ar[r]^{f_1}& X_{1}\ar[d]_{p_{1}}\ar@{-->}[ru]_{g_0'}      \\
L_{2}\ar[r]^{f_2}& X_{2}\ar[d]_{p_{2}}\ar@{-->}[rruu]_{g_1'}    \\
&\vdots\ar[d]_{p_{m-2}} \\
L_{m-1}\ar[r]^{f_{m-1}}& X_{m-1}\ar[d]_{p_{m-1}} \ar@{-->}[rrrruuuu]_{g_{m-2}'}   \\
& X_{m}.\ar@{-->}[rrrrruuuuu]_{g_{m-1}'}     \\
}
$$
Since $L_0$ is a wedge of spaces in $\F$ and $g_0\in Z_\F(X)$, we have $g_0\of f_0\simeq 0$
by Lemma \ref{lem:basic}(b).  Thus there is a map
$g_0': X_1\maprt{} X$ extending $g_0$.  
The same argument inductively 
defines  $g_i'$ for each $i$, and shows 
$g_{m-1} \cdots g_1 g_0\simeq g_{m-1}'\of ( p_{m-1}\cdots p_1 p_0 )$.
Now $g_{m-1} \cdots g_1 g_0\simeq 0$ since $X_m\equiv *$.  This proves the first assertion.

Next we let $m=\cl_\F(X)$, and show that $\kl_\F(X)\leq m$.  Let
$$
L_i\maprt{f_i} X_i\maprt{p_i} X_{i+1}
$$
for $0\leq i < m$ be an $\F$-cone length decomposition of $X$.  Set 
$$
\xymatrix{
h_i = (p_{m-1}p_{m-2} \cdots p_{i+1} )\of p_i : X_i\ar[rr] && X_m \equiv X\\ }
$$
for $i< m$ and $h_m = \mathrm{id}$.   Since $h_i=  h_{i+1}\of p_i$, Lemma \ref{lem:coflem}
yields cofiber sequences
$$
C_{p_i}\maprt{} C_{h_{i}} \maprt{} C_{h_{i+1}},
$$ 
for $0\leq i < m$.  This is a killing length decomposition of $X$.
To see this, observe that $C_{p_i}\equiv \s L_i$, which is a wedge of spaces in
$\F$ because $\F$ is closed under suspension.  
Furthermore,  $h_0: X_0\equiv *\maprt{} X$, so  $C_{h_0}\equiv X_m\equiv X$.
Finally,  $C_{h_m}\equiv 0$ because
$h_m=\mathrm{id}:X\maprt{} X$.  \halmos

\thm{Remarks}

\leftskip 28.5pt
\noindent\llap{(a)\qua}%
The notion of cone length has been extensively studied.
The version in Definition \ref{defn:killcone} is  similar to the one
given by Cornea in \cite{Cornea} (see (c) below).  It is precisely the same 
as the definition of ${\cal F}$-Cat given by Sheerer and Tanr\'e \cite{S-T}.
The $\F$-cone length $\cl_\F(X)$ can 
be regarded as the minimum number of steps needed to build the space $X$
up from a contractible space by attaching cones on wedges of spaces in $\F$.  
The notion of $\F$-killing length
is new and also appears in \cite{AMS} for the case $\F =\S$.  It can be regarded
as the minimum number of steps needed to destroy $X$ (i.e. go from $X$ to 
a contractible space) by attaching cones on wedges of spaces in $\F$.  
We note that Theorem \ref{thrm:nilbound} appears in 
\cite[Thm.$\,$3.4]{AMS}  
for the case $\F=\S$.
For the collection $\S$, it was shown in 
\cite[Ex.$\,$6.8]{AMS}  
that the inequalities in Theorem
\ref{thrm:nilbound} can be strict.

\noindent\llap{(b)\qua}
A space need not have a finite $\F$-killing length or $\F$-cone length decomposition.
For example, $\kl_\s(\cp^\infty)=\infty$ because all $2^n$-fold cup products vanish
in a space $X$ with $\kl_\s(X)\leq n$.  However, if
$X$ is a finite dimensional complex, then the process of attaching 
$i$-cells to the $(i-1)$-skeleton provides
$X$ with a $\S$-cone length decomposition.   Thus in this case, $\kl_\S(X)\leq \cl_\S(X)\leq \dim(X)$.
Since $\S\sseq\M\sseq \s$, it follows that $\kl_\s(X)\leq \kl_\M(X)\leq \kl_\S(X)$
and $\cl_\s(X)\leq \cl_\M(X)\leq \cl_\S(X)$, and so $\dim(X)$ is an upper bound
for all of these integers.  If $X$ is a $1$-connected finite dimensional complex, then 
a better upper bound for $\cl_\M(X)$ is the number of nontrivial positive-dimensional
integral homology groups of $X$.  This can be seen by taking a homology decomposition
of $X$ 
\cite[Chap.$\,$8]{Hi}.

\noindent\llap{(c)\qua}
It follows from work of Cornea \cite{Cornea} that the cone length of 
a space $X$, denoted $\cl(X)$, can be defined exactly like the $\s$-cone length $\cl_\s(X)$ above,
except that one does not require $L_0\in \s$.  It follows immediately that
$\cl(X)\leq \cl_\s(X)$.

\noindent\llap{(d)\qua}
The inequality $\kl_{\cal F}(X)\leq \cl_{\cal F}(X)$ also follows from work
of Sheerer and Tanr\'e since the function $\kl_{\cal F}$ satisfies
the axioms for ${\cal F}$-Cat 
\cite[Thm.$\,$2]{S-T}.

\leftskip0pt

We conclude this section by giving a few properties of killing length.

\bthm\label{thrm:killcof}
If $\F$ is any collection of spaces and $X\maprt{j} Y\maprt{q} Z$ is a cofiber
sequence, then 
$$
\kl_\F(Y)\leq \kl_\F(X) + \kl_\F(Z).
$$
\ethm

\pf Proof
Write $\kl_\F(X) = m$ and $\kl_\F (Z) = n$.  Let 
$$
L_i\maprt{f_i}X_i\maprt{} X_{i+1}
$$
for $0\leq i < m$ be a $\F$-killing length decomposition of $X$.
Set $g_0= j: X_0\equiv X \maprt{} Y$ and define $Y_1$ by the cofibration
$L_0\maprt{g_0 f_0} Y\maprt{} Y_1$.
By Lemma \ref{lem:coflem}, there is an auxilliary cofibration
$$
\xymatrix{
 C_{f_0} \ar[r]^{}\ar@{=}[d] &  C_{g_0 f_0} \ar[r]\ar@{=}[d] &  C_{g_0} \ar@{=}[d] \\
X_1\ar[r]^{g_1} & Y_1 \ar[r] & Z\\ }
$$
which defines $g_1$.
We proceed by induction: given $g_i:X_i\maprt{}Y_i$, let $Y_{i+1}$ be the
cofiber of the map ${g_i f_i} : L_0\maprt{}Y_i$
and use Lemma \ref{lem:coflem} to construct an auxilliary cofibration
$$
\xymatrix{
 C_{f_i} \ar[r]^{}\ar@{=}[d] &  C_{g_i f_i} \ar[r]\ar@{=}[d] &  C_{g_i} \ar@{=}[d] \\
X_{i+1}\ar[r]^{g_{i+1}} & Y_{i+1} \ar[r] & Z\\ }
$$
which defines $g_{i+1}$.  This defines cofiber sequences of the form
$L_j\to Y_j\to Y_{j+1}$ with $0\leq j < m$.
Since $X_m\equiv *$, the $(m+1)^{st}$ cofiber sequence, 
$X_m\to Y_m\to Z$,
shows that $Y_m\equiv Z$.  Now adjoin the $n$ cofiber sequences 
of a minimal $\F$-killing length decomposition of $Z$ to the first $m$
cofiber sequences to obtain an $\F$-killing length decomposition for $Y$
with length $n+m$. \halmos

Finally, we obtain an upper bound for $\kl_\s(X)$ and hence an upper bound for
$t_\s(X)$.  This provides a useful complement to Proposition \ref{prop:Et}
when $\cat(X)$ is not known.

\bprp\label{prop:killcat}
Let $X$ be an $N$-dimensional complex which is $(n-1)$-conn\-ected for some $n\geq 1$.
Then 
$$
\kl_\s(X)\leq \left\lceil \log_2\left( \frac{N+1}{n} \right) \right\rceil.
$$
\eprp

\pf Proof
We argue by induction on $\left\lceil \log_2\left( \frac{N+1}{n} \right)\right\rceil$. 
If $\left\lceil \log_2\left( \frac{N+1}{n} \right)\right\rceil=1$, then $N\leq 2n-1$.
It is well known that this implies that $X$ is a
suspension,  which means that $\kl_\s(X) =1$.  Now suppose 
$\left\lceil \log_2\left( \frac{N+1}{n} \right)\right\rceil= r$ and the result is 
known for all smaller values.  Let $X^k$ denote the $k$-skeleton of $X$, and consider
the cofiber sequence
$$
X^{2n-1} \maprt{} X \maprt{} X/X^{2n-1} .
$$
By Theorem \ref{thrm:killcof},
$\kl_\s(X) \leq \kl_\s(X^{2n-1}) + \kl_\s(X/X^{2n-1} )$.
The inductive hypothesis applies to $X^{2n-1}$  and 
to $X/X^{2n-1}$, so $\kl_\s(X) \leq 1 + (r-1) =r$.\halmos

\section{Distinguishing $Z_{{\cal F}}$ for different $\F$}

We have a chain of pointed sets
$$
Z_\Sigma(X,Y)\sseq Z_\M(X,Y)\sseq Z_\S(X,Y).
$$
Simple examples show that each of these containments
can be strict.  There are nontrivial phantom maps  $\s \cp^\infty\maprt{} S^4$ \cite{McG}.
These all lie in $Z_\M(\s \cp^\infty, S^4)$ because $\M\sseq {\cal P}$ (see Examples in Section 2), 
but not in $Z_\s(\s \cp^\infty, S^4)$, 
by Lemma \ref{lem:basic}(b).   For the other containment, the
Bockstein applied to the fundamental cohomology class of
$M(\ZZ/p, n )$ \cite{Baues} corresponds to a map $f:M(\ZZ/p, n) \maprt{} K(\ZZ/p, n+1)$.  If
$p$ is an odd prime, then $\pi_{n+1}(M(\ZZ/p, n))= 0$ 
\cite[pp.$\,$268--69]{Baues} 
so $f$ is in 
$Z_\S(M(\ZZ/p, n) , K(\ZZ/p, n+1))$. Since it is essential, $f$ cannot lie in 
$Z_\M(M(\ZZ/p, n) , K(\ZZ/p, n+1))$.

In these examples either the domain or the target is an infinite CW complex.  
Thus they leave open the possibility that if $X$ and  $Y$ are finite complexes,  
all of the pointed sets above are the same. We will give examples which show that, even for finite
complexes, these inclusions can be strict.   These examples are more
difficult to find and verify.  They are inspired by an example (due to Fred Cohen)
from \cite{F-M}.

Recall that if $p$ is an odd prime, then $S^{2n+1}_{(p)}$ is an H-space 
\cite{Adams}.
Moreover, if $f$ is in the abelian group $[\s^2 X,S^{2n+1}]$ then the order of 
$\lambda \of f\in [\s^2 X, S^{2n+1}_{(p)}]$ is either  infinite or a power of 
$p$.  

\blem\label{lem:moore}  Let $X$ be a finite complex and let $h:X\maprt{} 
S^{2n+1}$ be a map
such that for some odd prime $p$, $\lambda \of \s^2 h$ is nonzero and has finite order 
divisible $p$.  Then there is an $s>0$ such that the composite
$$
X \maprt{h} S^{2n+1} \maprt{i} M(\ZZ/p^s, 2n+1),
$$
where $i$ is the inclusion, is essential.
\elem

\pf Proof  Consider the diagram
$$
\xymatrix{
& S^{2n+1}\ar[d]_{p^s} \ar[r]^\lambda & S^{2n+1}_{(p)}\ar[d]^{p^s}\\
X\ar[r]^h \ar[rd]_{i\of h} & S^{2n+1}\ar[d]_{i} \ar[r]^\lambda & S^{2n+1}_{(p)}\ar[d]^j\\
&M(\ZZ/p^s,2n+1)\ar[r]^= &M(\ZZ/p^s,2n+1)\\ }
$$
in which the vertical sequences are cofibrations and $p^s$ denotes the map with degree $p^s$.    
If $i\of h\simeq 0$, then $j\of\lambda\of h=0$.  
It can be shown that $\s(\lambda\of h)$ lifts through the map
$p^s:S^{2n+2}_{(p)}\maprt{}S^{2n+2}_{(p)}$.   
Suspending once more, we obtain  a lift given by the dashed line in the diagram
$$
\xymatrix{
&& S^{2n+3}_{(p)}\ar[d]^{p^s}\\
\s^2X\ar[r]_{\s^2 h}\ar@{-->}[rru]& S^{2n+3}\ar[r]_\lambda & S^{2n+3}_{(p)}.\\ }
$$
Since $X$ is a finite complex, the torsion in $[\s^2 X,S^{2n+3}_{(p)}]$ is $p$-torsion
and has an exponent $e$. 
Since $S^{2n+3}_{(p)}$ is an H-space, the map $p^s$ induces multiplication by 
$p^s$
on $[\s^2 X, S^{2n+3}_{(p)}]$.
If $s\geq e$, then the image of $p^s:
[\s^2 X, S^{2n+3}_{(p)}]\maprt{}[\s^2 X, S^{2n+3}_{(p)}]$ cannot contain 
any nontrivial torsion.  But $\lambda \of\s^2 h$ 
is nonzero and has finite order.
Therefore the lift cannot exist, and so $i\of h\not\simeq 0$.
\halmos

\bthm\label{thrm:zeromap}
Let $X$ be a finite complex, let $p$ be an odd prime and let $g:X\maprt{} 
S^{2n+1}$
be an essential map.\ethm\sl

\leftskip28.5pt
\noindent\llap{\rm(a)\qua}%
Assume that $\pi_{2n+1}(X)$ is a finite group, and that 
$\lambda\of\s^2 g$ is nonzero with finite order divisible by $p^{n+1}$.
Then there is an $s>0$ such that the composite
$$
\xymatrix{
X\ar[r]_(.4){g} \ar@(ru, ul)[rrr]^{l} &
S^{2n+1}\ar[r]_{p^{n}}   &
S^{2n+1}\ar[r]_(.35){i}        &
M(\ZZ/p^s, 2n+1) \\ }
$$
is essential and $\pi_*(l)=0$.

\noindent\llap{\rm(b)\qua}%
Assume that $\pi_{k}(X)= 0$ for $k= 2n$ and  $2n+1$, and that
$\lambda \of \s^2 g$ is nonzero with finite order divisible by $p^{2n+1}$.
Then there is an $s>0$ such that the composite 
$$
\xymatrix{
X\ar[r]_(.4){g} \ar@(ru, ul)[rrr]^f &
S^{2n+1}\ar[r]_{p^{2n}}   &
S^{2n+1}\ar[r]_(.35){i}        &
M(\ZZ/p^s, 2n+1) \\ }
$$
is essential, and $\pi_*(f;G)= 0$ for any finitely generated
abelian group $G$.

\rm\leftskip0pt

\pf Proof
In part (a), the composition $\lambda\of \s^2 (p^{n}\of  g )$ 
has finite order divisible by $p$.  Therefore Lemma \ref{lem:moore} shows that 
$l=i\of p^{n}\of g$ is essential if $s$ is large enough.
Similarly, if $s$ is large enough, the map $f$ in part (b) is
essential.  From now on, we assume that $s$ has been so chosen.
We use the commutative diagram
$$
\xymatrix{
X \ar[r]^g\ar[rd]_{\lambda\of g} & S^{2n+1}\ar[d]^\lambda\ar[r]^{p^k} & 
S^{2n+1}\ar[r]^(.35){i}\ar[d]^\lambda & M(\ZZ/p^s, 2n+1)\ar[d]^=\\
& S^{2n+1}_{(p)} \ar[r]^{p^k} & S^{2n+1}_{(p)}\ar[r]^(.35){j}  & M(\ZZ/p^s, 2n+1).\\  }
$$
We take $k=n$ in part (a) and $k=2n$ in part (b).

\noindent{\bf Proof of} (a)\qua
Since $M(\ZZ/p^s, 2n+1)$ is $p$-local, there is only $p$-torsion to consider.
By results of Cohen, Moore and Neisendorfer
\cite[Cor.$\,$3.1]{CMN} 
the $p$-torsion in $\pi_*(S^{2n+1}_{(p)})$ has exponent
$n$.  Since $S^{2n+1}_{(p)}$ is an H-space,  
$p^n : S^{2n+1}_{(p)} \maprt{} S^{2n+1}_{(p)}$  annihilates all $p$-torsion in 
homotopy groups.
Thus $\pi_*(l)$ can be nonzero only in dimension $2n+1$.  But $\pi_{2n+1}(g)$
is a homomorphism from a finite group to $\ZZ$, so $\pi_*(l)=0$.

\noindent {\bf Proof of} (b)\qua  
It suffices to show that $\pi_m(f;G)=0$ for any cyclic group $G$;
by part (a) we need only consider $G=\ZZ/p^r$.
For each $r\geq 1$ and each $m\geq 0$, there is the
exact coefficient sequence 
\cite[Chap.$\,$5]{Hi} 
$$
0\maprt{} \Ext(  \ZZ/p^r, \pi_{m+1}(Y)  )
\maprt{} \pi_m(Y; \ZZ/p^r)
\maprt{}  \Hom (  \ZZ/p^r, \pi_{m}(Y)  ) \maprt{} 0.
$$
Let $Y=S^{2n+1}_{(p)}$. Since the $p$-torsion in $\pi_*(S^{2n+1}_{(p)})$ has
exponent $n$ \cite{CMN}, the exact sequence shows that the $p$-torsion in 
$\pi_m(S^{2n+1}_{(p)};\ZZ/p^r)$ has exponent at most
${2n}$ if $m\neq 2n$.  Thus the map $p^{2n} :  S^{2n+1}_{(p)}\maprt{}S^{2n+1}_{(p)}$
induces $0$ on the $m^{th}$ homotopy groups with coefficients in 
any finite abelian group if $m\neq 2n$.   
Taking $Y=X$ in the  coefficient  sequence, we have $\pi_{2n}(X;\ZZ/p^r)=0$.   
Therefore $\pi_*(f;G)=0$ for any finitely generated abelian group $G$. \halmos

We  apply this theorem to construct examples of finite
complexes which distinguish the various $Z_\F$.

Our first example shows that $Z_\M(X)$ can be different
from $Z_\S(X)$ even when $X$ is a finite complex.  Using the coefficient
exact sequence for homotopy groups, we find that 
$$
[M(\ZZ/p^r, 2n), S^{2n+1}] = \pi_{2n}(S^{2n+1};\ZZ/p^r)  \cong \ZZ/p^r
$$
for each $r$; this is a stable group.  Therefore, if $r>n$, there are 
essential maps $g:
M(\ZZ/p^r, 2n)\maprt{} S^{2n+1}$ with finite order divisible by $p^{n+1}$.  
Applying part (a) of Theorem \ref{thrm:zeromap}, we have the following example.

\thm{Example}
Let $r > n>1$.  For $p$ an odd prime and  $s$ large enough,  there are essential maps
$$
l:M (\ZZ/p^r, 2n)\maprt{} M(\ZZ/p^s, 2n + 1)
$$
such that $\pi_*(l)=0$.  Therefore by Lemma \ref{lem:self},
$$
Z_\S\bigl(   M(\ZZ/p^r, 2n) \wdg M(\ZZ/p^s, 2n + 1)     \bigr) \neq 0
$$
while, of course, 
$$
Z_\M\bigl(    M(\ZZ/p^r, 2n) \wdg M(\ZZ/p^s, 2n + 1)     \bigr)=0
$$
by Lemma \ref{lem:basic} (b).
It can be shown that any $s\geq r$ will
suffice in this example.

Freyd's generating
hypothesis \cite{Freyd} is the conjecture that no stably nontrivial map between finite complexes can
induce zero on stable homotopy groups.  The map $l$ in this example is stably nontrivial,
but our argument does {\it not} show that large suspensions of
$l$ induce zero on homotopy groups; the difficulty is that after two suspensions, 
$l$ factors through $p^n:S^{2n+3}_{(p)}\maprt{}S^{2n+3}_{(p)}$, which need not
annihilate all $p$-torsion.

Our second example is a map 
$f : \s^{2n-2}\left(\cp^{p^{2n+1}}/S^2\right) \maprt{} M(\ZZ/p^s, {2n+1})$ 
which we use to show that $Z_\M(X)$
can be different from $Z_\Sigma(X)$ when $X$ is a finite complex.  We need some
preliminary results to show that Theorem \ref{thrm:zeromap} applies to this situation.

\blem\label{lem:mcg}
Let $f:\s^{n+1}\cp^m\maprt{} S^{n+3}$.  The degree of 
$f|_{\s^{n+1}S^2}$ is divisible by 
$\lcm(1,\ldots, m)$, the least common 
multiple of $1,\ldots, m$.
\elem

\pf Proof
We may assume that $f$ is in the stable range.
If $f|_{\s^{n+1}S^2}$ has degree $d$, then
$$
\s^{n+1} \cp^m \maprt{f} S^{n+3}\inclds \s^{n+1} \cp^m.
$$
has degree $d$ in $H_{n+3}(\s^{n+1} \cp^m) $ and is trivial in all other dimensions.
According to McGibbon 
\cite[Thm.$\,$3.4]{McG},
$d$ is divisible by $\lcm(1,\ldots, m)$.\halmos

\bprp\label{prop:apply}
The image of the $n$-fold suspension map
$$
\s^n: [ \cp^{p^t}/S^2 , S^{3}]
\maprt{} [\s^{n}( \cp^{p^t}/S^2 ), S^{n+3}]
$$
contains elements of order $p^t$ for every $n\geq 1$ and $t\geq 1$.
\eprp

\pf Proof
Write $m=p^t$ and examine the commutative diagram
$$
\xymatrix{
[\s\cp^{m}, S^3]\ar[r]\ar[d]^{\s^n} &[\s S^2, S^3]\ar[r]\ar[d]^{\s^n} &[\cp^{m} /S^2 , S^3]\ar[d]^{\s^n}\\
[\s^{n+1}\cp^{m}, S^{n+3}]\ar[r]\ar[d]^{\lambda_*} &
      [\s^{n+1} S^2, S^{n+3}]\ar[r]\ar[d]^{\lambda_*} &[\s^n (\cp^m /S^2) , S^{n+3}]\ar[d]^{\lambda_*}\\ 
[\s^{n+1}\cp^{m}, S^{n+3}_{(p)}]\ar[r] &
    [\s^{n+1} S^2, S^{n+3}_{(p)}]\ar[r]& 
    [\s^n (\cp^{m} /S^2) , S^{n+3}_{(p)} ].\\  }
$$
To show that the image of 
$\lambda_*\of \s^n: [\cp^{m} /S^2 , S^3]\maprt{}[\s^n (\cp^{m} /S^2) , S^{n+3}_{(p)}]$
contains elements of order $p^t$, we modify the above diagram
as follows:  the image and cokernel of $[\s\cp^m,S^3]\maprt{} [\s S^2, S^3]\cong \ZZ$
are $k\ZZ$ and $\ZZ/k$, respectively, for some integer $k$; similarly
for $[\s^{n+1}\cp^m,S^{n+3}]\maprt{} [\s^{n+1} S^2, S^{n+3}]$
and $[\s\cp^m,S^3_{(p)}]\maprt{} [\s S^2, S^3_{(p)}]$.  Thus we have a commutative
diagram with exact rows
%
%
$$
\xymatrix{
k\ZZ \ar[r]\ar[d] &\ZZ \ar[r]\ar[d]^\cong &\ZZ/ k\ar[d]\\
l\ZZ \ar[r]\ar[d]^\lambda &\ZZ\ar[r]\ar[d]^\lambda &\ZZ/l \ar[d]^\lambda\\
l_{p}\ZZ_{(p)} \ar[r] &\ZZ_{(p)} \ar[r]& \ZZ/l_{p} \\  }
$$
for some integers $k$, $l$ and $l_{p}$, where $l_p$ is the largest power of $p$
which divides $l$.  Lemma \ref{lem:mcg} shows that  $l_p$ is divisible by $p^t$.  
The composite
$\ZZ/k\maprt{} \ZZ/l\maprt{} \ZZ/l_p$ is surjective, and this completes the proof.\halmos

It follows from Proposition \ref{prop:apply} that part (b) of 
Theorem \ref{thrm:zeromap} applies to the space 
$\s^{2n-2} ( \cp^{p^{2n+1}}/S^2 )$ for each $n>1$, and so we obtain our second example.

\thm{Example}
For each odd prime $p$ and each $n\geq 1$, there is an 
$s>0$ such that there are essential maps
$$
f: \s^{2n-2} \left( \cp^{p^{2n+1}}/S^2\right) \maprt{} M(\ZZ/p^s , 2n+1)
$$
which induce zero on homotopy groups with coefficients.  Therefore,
$$
Z_\M  \left(\s^{2n-2} \left(\cp^{p^{2n+1}}/S^2\right)\wdg M\left(\ZZ/p^s , 2n+1\right) \right)\neq 0
$$
while, of course, 
$$
Z_\s  \left(\s^{2n-2} \left(\cp^{p^{2n+1}}/S^2\right)\wdg M\left(\ZZ/p^s , 2n+1\right)
\right)=0.
$$
The map $f$ can be chosen to be
stably nontrivial. As in the previous example, the suspensions of $f$ might not be trivial on 
homotopy groups with coefficients.

Finally, let $A = \s^{2n-2} (\cp^{p^{2n+1}}/S^2)$,
$B = M(\ZZ/p^r, 2n)$ and $C= M(\ZZ/p^s, 2n + 1)$ for $s$ large.
Then 
$$
Z_\s (A\wdg B\wdg C)  < Z_\M (A\wdg B\wdg C) < Z_\S(A\wdg B\wdg C),
$$
so both of these inequalities can be strict for a single finite complex.

\section{Projective spaces}

We show that for projective spaces $\fp^n$ with $\FF = \RR,\CC$ or $\HH$, 
$$
Z_\s(\fp^n) = Z_\M(\fp^n)=Z_\S(\fp^n),
$$
and we completely determine these sets for $\FF=\RR$ and $\CC$ and all $n$.  We also
determine $t_\S(\hp^n)$, for $n\leq 4$.

\subsection{General facts}
We first prove some 
general results that will be applied later.

\bprp\label{prop:ziss}
If $\s\om X\equiv \bigvee S^{n_\alpha}$,
then  $Z_\S(X,Y) = Z_\M(X,Y)=Z_\s(X,Y)$ for any space $Y$.
\eprp

\pf Proof
Let $f\in Z_\S(X,Y)$.\ The map
$\om f$ is adjoint to the composition
$\s\om X\! \maprt{\nu}$\break$ X \maprt{f} Y$.
Since $\s\om X\equiv \bigvee S^{n_\alpha}$, $f\of \nu\simeq 0$, and
so $\om f\simeq 0$.  Thus $f\in Z_\s(X,Y)$.\halmos

By Lemma \ref{lem:E}, the condition $Z_\S(X,Y)=Z_\s(X,Y)$ is equivalent to
the condition that if $f:X\maprt{}Y$ induces zero on homotopy groups, then $\om f=0$.

Proposition \ref{prop:ziss} applies to $X =S^{n+1}$  because,
by James \cite{J}, $\s \om S^{n+1}\equiv \bigvee_{k=1}^\infty S^{nk+1}$.
Of course $Z_\S(S^{n+1},Y)=0$.  Since 
$\s(A\cross B)\equiv \s A\wdg \s B\wdg \s(A\smsh B)$ for any $A$ and 
$B$ 
\cite[11.10]{Hi},  
James's result
allows us to apply Proposition \ref{prop:ziss} to any space
whose loop space splits as a finite type product of spheres and loop spaces on spheres.  
Moreover, if $X$ and $X'$ both satisfy the hypotheses of Proposition \ref{prop:ziss},
then so does $X\cross X'$.

For $\FF=\RR, \CC$ or  $\HH$, let $d= 1, 2$ or $4$,  respectively.  
For each $n\geq 1$  there is a homotopy equivalence
$\om \fp^n \equiv  S^{d-1}\cross \om S^{(n+1)d -1}$.  This is 
a direct consequence of \cite[Thm.$\,$5.2]{E-H}, which applies
even in the case $d=1$.

\bcor
For $\FF=\RR,\CC$ or $\HH$ and each $n\geq 1$, 
$Z_\S( \fp^n ) =Z_\M(\fp^n)= Z_\s(\fp^n)$.
\ecor

Another corollary of Proposition \ref{prop:ziss} applies to intermediate wedges of
spheres.
For spaces $X_1,X_2, \ldots , X_n$, the elements $(x_1,\ldots ,x_k)\in X_1\cross \cdots \cross X_k$
with at least $j$ coordinates equal to the base point form a subspace
$T_j(X_1, \ldots, X_k)\sseq X_1\cross \cdots \cross X_k$.
Porter has shown 
\cite[Thm.$\,$2]{Porter}
that
$\om T_j(S^{n_1}, \ldots, S^{n_k})$ has the homotopy type of a product of 
loop spaces of spheres
for each $0\leq j\leq k$.  Our previous discussion establishes
the following.

\bcor  \label{cor:spheres}
For any $n_1, \ldots, n_k \geq 1$ and any $0\leq j\leq k$,  
$$
Z_\S( T_j( S^{n_1}, \ldots,  S^{n_k})) =Z_\M( T_j( S^{n_1}, \ldots,  S^{n_k}))=
Z_\s( T_j( S^{n_1}, \ldots,  S^{n_k}) ).
$$
\ecor

\thm{Remarks}
\begin{enumerate}
\item[(a)]   
Taking $j=0$ in Corollary \ref{cor:spheres}, we deduce from Corollary \ref{cor:spheret}
that
$$
t_\S (S^{n_1}\cross \cdots \cross S^{n_k}) \leq \big\lceil \log_2\left( k+1\right) \big\rceil .
$$
This reproves 
\cite[Prop.$\,$6.2]{AMS}
by a different method.
\item[(b)]
It is proved in 
\cite[Prop.$\,$6.5]{AMS}
that for any positive integer
$n$, there is a finite product of spheres $X$ with $t_\S(X)=n$.  By Corollary \ref{cor:spheres},
the same is true for $t_\s(X)$ and $t_\M(X)$.  Thus the integers $t_\F(X)$ for 
$\F=\S, \M$ or $\s$ and any $X$ take on all positive integer values.
\end{enumerate}

Finally, we observe that the splitting of $\om \fp^n$ gives a useful criterion for deciding when a map
$f:\fp^n\maprt{}Y$ lies in $Z_\S(\fp^n, Y)$.

\bprp\label{prop:criterion}
Let $i$ be the inclusion $S^d=\fp^1\inclds \fp^n$, and 
let $p:S^{(n+1)d-1}\maprt{} \fp^n$ be the Hopf fiber map.  Then the map
$$(i,p) : S^d \wdg S^{(n+1)d-1}\maprt{} \fp^n$$
induces a surjection on homotopy groups.  Therefore, a map
$f: \fp^n\maprt{}Y$ satisfies $\pi_*(f)=0$ if and only if $f\of i=0$
and $f\of p=0$.
\eprp

\subsection{Complex projective spaces}
 
Next  we show that certain skeleta $X$ of Eilenberg-MacLane spaces
have the property that $Z_\S(X) = 0$.  
We apply this to $\cp^n$ and  $\s^n \cp^2$ for each $n$.

Let $G$ be a finitely generated abelian group.
Give the Eilenberg-MacLane space $K(G,n)$ with $n\geq 2$  a 
homology decomposition 
\cite[Chap.$\,$8]{Hi} 
and denote the $m^{th}$ section by $K(G,n)_m$.  Thus
$K(G,n)$ is filtered 
$$
*\sseq K(G,n)_n\sseq \cdots\sseq K(G,n)_m\sseq \cdots \sseq K(G,n)
$$
and there are cofiber sequences
$$M( H_{m+1}(K(G,n)), m)\maprt{} K(G,n)_m \maprt{} K(G,n)_{m+1}.$$

\bthm\label{thrm:Km}
If the group $H_{m}(K(G,n))$ is torsion free and 
$H_{m+1}(K(G,n))=0$, then $Z_\S(K(G,n)_m) = 0$.
\ethm

\pf Proof
We write $X= K(G,n)_m$.  Then $H_k(K(G,n),X) = 0$ for $k\leq m+1$.  By 
Whitehead's theorem 
\cite[Thm.$\,$7.13]{Wh}, 
the induced map
$\pi_k(X) \maprt{} \pi_k( K(G,n))$ is an isomorphism for $k\leq m$.  
Since $H_{m}(K(G,n))$ is torsion free,
$X$ has dimension at most $m$, and so $X$ has a CW decomposition 
$$
\bigvee S^n = X^n \sseq X^{n+1}\sseq \cdots \sseq X^{m} \equiv X.
$$
For $f\in Z_\S(X)$ we prove
by induction on $k$ that $f$ factors through $X/X_k$ for each $k\leq m$.    The first step is
trivial  since $\pi_n(f) = 0$ implies $f|_{X_n}\simeq 0$.  Inductively,
assume that $f$ factors through $X/X_k$ with $n\leq k < m$.  There is
a cofibration
$$
\bigvee S^{k+1}\equiv X_{k+1}/X_k\maprt{} X/X_{k} \maprt{} X/X_{k+1}.
$$
Since $n< k+1\leq m$, it follows that $\pi_{k+1}(X)\cong \pi_{k+1}(K(G,n)) = 0$,
so $f$ extends to $X/X_{k+1}$.
Taking $k=m$, we find $f\simeq 0$.  \halmos

\thm{Remark}
Clearly,  $\pi_k(K(G,n)_m)=0$ for $n < k < m$.  The hypotheses
in Theorem \ref{thrm:Km} are  needed to conclude further 
that $\pi_m(K(G,n)_m)=0$.

\medskip

As an application of Theorem \ref{thrm:Km}, we have the following calculations.

\bthm  \hfill
\begin{enumerate}
\item[(a)]  $Z_\F(\cp^n) = 0$ for each $n\geq 1$ and each $\F=\s,\M$ or $\S$.
\item[(b)]  $Z_\F(\Sigma^{n} \cp^2) = 0$ for each $n\geq 1$ and each $\F=\s,\M$ or $\S$.
\end{enumerate}
\ethm

\pf Proof
By Proposition \ref{prop:ziss} it suffices to consider the case $\F=\S$.
Since $\cp^\infty = K(\ZZ, 2)$ and the $\cp^n$ are the sections of a homology
decomposition of $\cp^\infty$,
part (a) follows from Theorem \ref{thrm:Km}.  Recall from \cite{E-M} that for $n\geq 2$ 
$$
H_k(K(\ZZ,n))  =     \bigg\{ 
\begin{array}[]{ll}
\ZZ      & \mbox{if $k=n$ or $n+2$} \\
0        & \mbox{if $k=n+1$ or $n+3$.}
\end{array}                                    
$$
Since $\mathrm{Sq}^2$ is nontrivial on $H^n(K(\ZZ,n);\ZZ/2)$,
we have $K(\ZZ,n)_{n+2}\equiv \s^{n-2}\cp^2$. 
Thus Theorem \ref{thrm:Km} applies to $\s^{n-2}\cp^2$.
\halmos

This theorem immediately shows that
$t_\F(\cp^n) = t_\F(\Sigma^{n} \cp^2) = 1$ for $\F=\s,\M$ or $\S$ and each $n\geq 1$.

\subsection{Real projective spaces}

In this subsection we completely  calculate  $Z_\S(\rp^n)$.  
By the Hopf-Whitney theorem
\cite[Cor.$\,$6.19]{Wh} 
$[\rp^{2n}, S^{2n}]\cong
H^{2n}(\rp^{2n})\cong\ZZ/2$. The unique nontrivial map $q:\rp^{2n}\maprt{}S^{2n}$ is 
the quotient map obtained
by factoring out $\rp^{2n-1}$.  Let $f_{2n}$ denote the composite
$\rp^{2n}\maprt{q}S^{2n}\maprt{p}\rp^{2n}$
where $p$ is the universal covering map.

\bthm  \label{thrm:rpn}
For $\F=\s, \M$ or $\S$ and each $n\geq 1$, 
\begin{enumerate}
\item[\rm(a)]  $Z_\F(\rp^{2n-1}) = 0$
\item[\rm(b)]  $Z_\F(\rp^{2n}) = \{ 0, f_{2n}\}$.
\end{enumerate}
\ethm

\pf Proof
Let $f:\rp^{m}\maprt{}\rp^m$ with $\pi_1(f) =0$.  
Because $\pi_k(\rp^{m})=0$
for $1<k< m$, an argument similar to the proof of Theorem \ref{thrm:Km}
shows that $f$ must factor through $q:\rp^{m}\maprt{} S^{m}$.  For $m>1$,
any map $S^m\maprt{}\rp^m$ lifts through $p:S^m\maprt{} \rp^m$.  Thus there is
a map $g:S^m\maprt{}S^m$ of degree $d$ which makes the following diagram commute
$$
\xymatrix{
S^{m}\ar[r]^p \ar[rd]_{q\of p} & \rp^{m}\ar[r]^f\ar[d]^q & \rp^{m}\\
& S^{m} \ar[r]^g & S^{m}.\ar[u]_p\\ }
$$

First let $m=2n-1$.  We may assume $n>1$.   The composite $q\of p :
S^{2n-1}\maprt{}S^{2n-1}$  is known to have degree $2$. Since $\pi_i(p)$ is an isomorphism for $i>1$, 
$f\of p$ represents 
$2d\in\ZZ\cong \pi_{2n-1}(\rp^{2n-1})$.  If $f\in Z_\S(\rp^{2n-1})$, 
then $d$ must be $0$, and so $f=0$.  This proves (a).

Now take $m= 2n$.  The composite $q\of p:S^{2n}\maprt{}S^{2n}$ is trivial because it
is zero on homology. 
Therefore $f\of p=0$, and since $\pi_1(f)=0$, Proposition \ref{prop:criterion}
shows that $f\in Z_\S(\rp^{2n})$.
Since $\rp^{2n}$ is connected, there is a bijection
$$
\xymatrix{
p_* : [\rp^{2n}, S^{2n}]\ar[rr]^(.3){\cong} && 
\big\{ f \| f\in [\rp^{2n}, \rp^{2n}] ,\, \pi_1(f)= 0\big\} = Z_\S(\rp^{2n}).\\  }
$$
Since $[\rp^{2n},S^{2n}] = \{ 0, q\}$ as noted above, 
$Z_\S(\rp^{2n}) = \{ 0, f_{2n}\}$, where $f_{2n}= p\of q$. \halmos

\thm{Remark}  This argument actually shows that, if $\pi_k(Y)=0$ for $1< k< 2n+1$,
there is a bijection between
$Z_\S(\rp^{2n+1},Y )$ and the set of elements $\alpha\in \pi_{2n+1}(Y)$ such that $2\alpha =0$.

\bcor  
For each $n\geq 1$, 
\begin{enumerate}
\item[\rm(a)]  $t_\F(\rp^{2n-1}) = 1$ for $\F= \s,\M$ and $\S$.
\item[\rm(b)]  $t_\F(\rp^{2n}) = 2$ for $\F= \s,\M$ and $\S$.
\end{enumerate}
\ecor

\pf Proof
It suffices to prove part (b) for $\F=\S$.
Since $Z_\S(\rp^{2n})\neq 0$, $t_\S(\rp^{2n})$\break$\geq 2$.  The only possibly nonzero
product in this semigroup is $f_{2n}\of f_{2n}$.  But this is zero because 
$Z_\S(\rp^{2n})$ is nilpotent by Theorem \ref{thrm:nilbound}.\halmos

\subsection{Quaternionic projective spaces}

The quaternionic projective spaces are not  skeleta of  Eilenberg-MacLane spaces,
and it is much more difficult to compute their nilpotency.

Let $f\in [\hp^{n+1}, \hp^{n+1}]$, and assume that $f$ is cellular.
Then $f|_{\hp^n}: \hp^n\maprt{}\hp^n$ and the homotopy
class  $f|_{\hp^n}$ is well defined. 

\blem\label{lem:hprestrict}
If $f\in Z_S(\hp^{n+1})$, then $f|_{\hp^n} \in Z_\S(\hp^n)$.
\elem

\pf Proof
Let $f\in Z_S(\hp^{n+1})$ and let $g= f|_{\hp^{n}}$.  Consider the diagram
$$
\xymatrix{
& S^{4n+3} \ar@{-->}[rr]^h\ar[d]_p && S^{4n+3}\ar[d]^p\\
S^4\ar[r]^i & \hp^n\ar[d]_j\ar[rr]^g && \hp^n\ar[d]_j\ar[rd]^l\\
&\hp^{n+1}\ar[d]_q \ar[rr]^f &&\hp^{n+1}\ar[r]^m\ar[d]^q & \hp^\infty\\
&S^{4n+4}\ar@{-->}[rr]^{\s h} &&S^{4n+4}.\\ }
$$
where $i, j, m$ and $l$ are inclusions.
Since $S^{4n+3}\maprt{p} \hp^n\maprt{l}\hp^\infty$
can be regarded as a fibration and $l\of (g\of p) = m\of ( f\of ( j\of p ))= 0$, 
it follows that $g\of p$
lifts to the map $h$.  Since $f\in Z_\S(\hp^{n+1})$, $f$ induces 
zero on $H^4(\hp^{n+1})$, and hence is zero in cohomology.  
Therefore $\s h$ is zero in cohomology and hence is trivial.  Thus $h=0$,
so $g\of p=0$.  Also, $g\of i=0$, so $g\in Z_\S(\hp^n)$ by Proposition \ref{prop:criterion}. \halmos

Next we indicate how we will apply Lemma \ref{lem:hprestrict}.
If $Z_\S(\hp^n) = 0$ and  $f\in Z_\S(\hp^{n+1})$, then 
$f|_{\hp^n} = 0$, so $f$ factors through $q:\hp^{n+1}\maprt{}S^{4n+4}$.  
By Proposition \ref{prop:criterion}, if $i:S^4\inclds \hp^{n+1}$, 
then  $\pi_{4n+4}(i)$ is surjective, so $f$ factors as in the diagram
$$
\xymatrix{
\hp^{n+1} \ar[r]^f \ar[d]_q & \hp^{n+1} \ar[r]^m & \hp^\infty\\
S^{4n+4}\ar[r]^{g}& S^4.\ar[u]_i\ar[ru]_l \\  }
$$
By cellular approximation, $f$ is essential if and only if $m\of f$ is
essential.  The map $l\of g : S^{4n+4}\maprt{}\hp^{\infty}$ is adjoint
to a map $g':S^{4n+3}\maprt{}S^3$, which in turn is adjoint to
$l\of \s g'$.  By cellular approximation again,  
$i\of g = i\of\s g'$, so
we may  assume that $g$ is in the image of the 
suspension $\s:\pi_{4n+3}(S^3)\maprt{} \pi_{4n+4}(S^4)$.

The proof of our main result about quaternionic projective spaces
requires some detailed information about homotopy groups of spheres.
Since we refer to Toda's book \cite{Toda} for this information, we 
use his notation here.   For example, $\eta_k: S^{k+1}\maprt{} S^k$ and 
$\nu_k: S^{k+3}\maprt{} S^k$ are suspensions of the Hopf fiber maps.

\bthm\label{thrm:hpn}\hfill
\begin{enumerate}
\item[\rm(a)]  $Z_\F(\hp^n)= 0$ for $\F=\S,\M$ or $\s$ and  $n=1,2$ and $3$
\item[\rm(b)]  $Z_\F(\hp^4)\neq 0$ for $\F=\S,\M$ or $\s$.
\end{enumerate}
\ethm

\pf Proof
First $\hp^1=S^4$, so $Z_\S(\hp^1)=0$. 
If $f\in Z_\S(\hp^2)$, then there is a commutative diagram
$$
\xymatrix{
\hp^{2}\ar[r]^f\ar[d]_q     &   \hp^{2} \ar[r]^m      & \hp^\infty\\
S^{8}\ar[r]^{g}\ar[d]_{\nu_5}  &   S^4\ar[u]_i \ar[ru]_l \\ 
S^5\ar[ru]_{\eta_4} \\ }
$$
in which the vertical sequence is a cofibration.
If $g=0$, then $f=0$, so we may assume that $g\neq 0$.
We know that $\pi_8(\hp^\infty )  \cong \pi_7(S^3)\cong \ZZ/2$,
generated by $\eta_3\of \nu_4$ 
\cite[p.$\,$43--44]{Toda}.  
Thus we can take $g= \eta_4\of \nu_5$.  Since $\nu_5\of q =0$, we conclude that $g\of q=0$,
so $f=0$.  This shows that $Z_\F(\hp^2)=0$.

The proof that $Z_\S(\hp^3)=0$ is similar. Let $f\in Z_\S(\hp^3)$ and  apply Lemma
\ref{lem:hprestrict} to get a similar factorization.  The resulting map 
$g:S^{12}\maprt{}S^4$ is either $\s \epsilon_3$  or $0$ 
\cite[Thm.$\,$7.1]{Toda}.   
If $g=\s\epsilon_3$, then
results of 
\cite[(2.20a)]{James2} 
and 
\cite[Thm.$\,$7.4]{Toda} 
show that $f\of p\neq 0$.  Thus
$g=0$ and so $f=0$.

For part (b), we make use of the diagram preceding Theorem \ref{thrm:hpn}  and the fact that
$g$ can be taken to be a suspension map.  If $f\in Z_\S(\hp^4)$, then 
we have
$$
\xymatrix{
S^{19}\ar[r]^p\ar[rd] &\hp^{4}\ar[r]^f\ar[d]_q     &   \hp^{4} \ar[r]^m      & \hp^\infty\\
&S^{16}\ar[r]^{g}  &   S^4.\ar[u]_i \ar[ru]_l \\  }
$$
According to Toda \cite{Toda}, $\pi_{16}(S^4)= \s( \pi_{15}(S^3))\cong \ZZ/2 \oplus\ZZ/2$.
By 
\cite[(2.20a)]{James2}, 
$q\of p$ is $4\nu_{16}\in \pi_{19}(S^{16})$.
Then $g\of (4\nu_{16})= 4g \of \nu_{16}$ because $\nu_{16}$ is a suspension.  Since $4g =0$, any map
$\hp^4\maprt{}\hp^4$ which factors through $q$ lies in $Z_\S(\hp^4)$.   Marcum and Randall
show in \cite{M-R} that the map
$$
\xymatrix{
((i \of \s\nu')\of \mu_7)\of q  : \hp^4 \ar[rr] && \hp^4\\ }
$$
is essential, where $\nu'\in \pi_6(S^3)$ generates the $2$-torsion and 
$\mu_7\in \pi_{16}(S^7)$ generates a $\ZZ/2$ summand 
\cite[Thm.$\,$7.2]{Toda}.  
Thus  $Z_\S(\hp^4)\neq 0$, and so
$t_\S(\hp^4)\geq 2$. 
\halmos

As before, we obtain the  nilpotency.

\bcor\label{cor:thpn}\hfill
\begin{enumerate}
\item[\rm(a)]  $t_\F(\hp^1) =  t_\F(\hp^2) = t_\F(\hp^3) = 1$ for $\F=\S,\M$ or $\s$
\item[\rm(b)]  $t_\F(\hp^4) = 2$ for $\F=\S,\M$ or $\s$.
\end{enumerate}
\ecor

\pf Proof
It suffices to prove that $t_\S(\hp^4)\leq 2$.  Suppose $f,g\in Z_\S(\hp^n)$.
The proof of Theorem \ref{thrm:hpn} shows that $f$ factors through $S^{16}$.
Now $g\of f=0$ because $g\in Z_\S(\hp^4)$.\halmos

\section{H-spaces}

In this section we study the nilpotency of H-spaces $Y$ mod $\s$.   We make
calculations for specific Lie groups such as $SU(n)$ and $Sp(n)$ and show 
that $Z_\s$ is nontrivial in these cases.  
If $Y$ is an H-space, the Samelson product of $\alpha\in \pi_m(Y)$
and $\beta\in \pi_n(Y)$ is written $\< \alpha, \beta\> \in \pi_{n+m}(Y)$
\cite[Chap.$\,$X]{Wh}.

We first give a few general results which are needed later.

\blem\label{lem:5.3}
If $Y$ is an H-space and $\<\alpha,\beta\>\neq 0$ for some  
$\alpha\in \pi_m(Y)$ and $\beta\in \pi_n(Y)$, then $[S^m\cross S^n, Y]$
is not abelian.
\elem

\pf Proof
The quotient map $q:S^m\cross S^n\maprt{} S^m\smsh S^n\equiv S^{m+n}$
induces a monomorphism 
$q^*:[S^m\smsh S^n, Y]\maprt{} [S^m\cross S^n, Y]$
such that $q^*\< \alpha,\beta\>  = [\alpha p_1, \beta p_2]$,
the commutator of $\alpha p_1$ and $\alpha p_2$.\halmos

It is well known that if an H-space $Y$ is a finite complex, then it has the
same rational homotopy type as a product of spheres
$S^{2n_1 -1}\cross \cdots\cross S^{2n_r -1}$ with $n_1\leq \cdots\leq n_r$.
If $p$ is an odd prime such that
$$
Y_{(p)} \equiv \left( S^{2n_1 -1}\cross \cdots\cross S^{2n_r -1}\right)_{(p)} \equiv
S^{2n_1 -1}_{(p)}\cross \cdots\cross S^{2n_r -1}_{(p)},
$$
then $p$ is called a \term{regular prime} for $Y$.  If $Y$ is a simply-connected compact Lie
group, then $p$ is regular for $Y$ if and only if $p\geq n_r$ 
\cite[Sec.$\,$9-2]{Kane}.

We need a second
product decomposition for $p$-localized Lie groups.  By 
\cite[Sec.$\,$2]{Mimura}  
there are fibrations
$S^{2k+1}\maprt{} B_k(p) \maprt{} S^{2k + 2p -1}$
for $k= 1, 2, \ldots$. An odd prime $p$
is called \term{quasi-regular} for the H-space $Y$ if 
$$
Y_{(p)}\equiv \left( \prod_i S^{2n_i-1}\cross \prod_j B_{m_j}(p)\right)_{(p)}.
$$

\subsection{The groups $SU(n)$ and $Sp(n)$}
We apply the notions of regular and quasi-regular primes to the Lie group $SU(n)$,
which has the rational homotopy type of $S^3\cross S^5 \cross\cdots \cross S^{2n-1}$,
and to the Lie group $Sp(n)$, which has the rational homotopy type of $S^3\cross S^7\cross\cdots
\cross S^{4n-1}$.  
It is well known 
\cite[Thm.$\,$4.2]{Mimura}  
that if $p$ is an odd prime then 
\begin{enumerate}
\item[(a)] $p$ is regular for $SU(n)$ if and only if $p\geq n$; $p$ is
quasi-regular for $SU(n)$ if and only if $p> \frac{n}{2}$
\item[(b)]  $p$ is regular for $Sp(n)$ if and only if $p\geq 2n$; $p$ is
quasi-regular for $Sp(n)$ if and only if $p>n$.
\end{enumerate}
It is also known 
\cite[Thm.$\,$1]{Bott} 
that if $n\geq 3$
and  $r+s +1 = n$, there are generators  
$\alpha\in \pi_{2r+1}(SU(n))\cong \ZZ$, $\beta\in \pi_{2s+1}(SU(n))\cong\ZZ$ and 
$\gamma\in \pi_{2n}(SU(n))\cong\ZZ/n!$ such that  
$\< \alpha, \beta \> = r! s! \gamma $.
If $p$ is an odd prime and  
$\alpha'\in \pi_{2r+1}(SU(n)_{(p)})$, $\beta'\in \pi_{2s+1}(SU(n)_{(p)})$
and $\gamma'\in \pi_{2n}(SU(n)_{(p)})$
are the images of $\alpha$, $\beta$  and $\gamma$ under the localization homomorphism
$\lambda_*:\pi_*(SU(n)) \maprt{} \pi_{*}(SU(n)_{(p)}) \cong \pi_{*}(SU(n))_{(p)}$,
then
$$
\< \alpha' , \beta'\> = r! s! \gamma'\in \pi_{2n}(SU(n))_{(p)}\cong \ZZ/n! \otimes \ZZ_{(p)}.
$$
Now we  prove the main result of this section.

\bthm  The groups 
\begin{enumerate}
\item[\rm(a)] $[SU(n),SU(n)]$ for $n\geq 5$ and
\item[\rm(b)] $[Sp(n),Sp(n)]$ for $n\geq 2$ 
\end{enumerate}
are not abelian.
\ethm

\pf Proof
Consider $SU(n)$ for $n\geq 5$ and let $p$ be the
largest prime such that $\frac{n}{2} < p < n$.  If $n\geq 12$, then it follows from Bertrand's
postulate 
\cite[p.$\,$137]{Bertrand}
that there are two  primes $p$ and $q$ 
such that $\frac{n}{2}<q < p< n$. This implies that $2n+6 < 4p$.  For $5\leq n < 12$, and 
$n\neq 5, 7, 11$, it is easily verified that $2n+6 < 4p$.

Assume that $n\geq 5$ and that $n\neq 5,7$ or $11$.
Since $p>\frac{n}{2}$, it follows that $p$ is quasi-regular for $SU(n)$.
Since $2n+6 < 4p$,  the spheres $S^{2n-2p+3}$ and $S^{2p-3}$ both appear in the 
resulting product decomposition. Thus we have
$$
SU(n)_{(p)} \equiv \left( B_1(p) \cross \cdots\cross  B_{n-p}(p)
\cross S^{2n-2p+3} \cross \cdots \cross S^{2p-3}\cross S^{2p-1}\right)_{(p)}.
$$
Assume $[SU(n),SU(n)]$ is abelian. Then $[SU(n)_{(p)},SU(n)_{(p)}]$ is abelian, and 
therefore $[ S^{2n-2p+3}\cross S^{2p-3}, SU(n)_{(p)}]$ is abelian.

There are 
$\alpha'\in \pi_{2n-2p+3}(SU(n)_{(p)})$, $\beta'\in\pi_{2p-3}(SU(n)_{(p)})$ 
and\hfil\break 
$\gamma'\in \pi_{2n}(SU(n)_{(p)})$  so that
$$
\< \alpha', \beta' \> = (n-p +1)! (p-2)!\gamma'
$$
in $\pi_{2n}(SU(n)_{(p)})\cong \ZZ/n!\otimes \ZZ_{(p)}\cong \ZZ/p$.
Since $\gamma'$ is a generator of $\ZZ/p$, we have $\<\alpha',\beta'\>\neq 0$.  By 
Lemma \ref{lem:5.3}, $[ S^{2n-2p+3}\cross S^{2p-3}, SU(n)_{(p)}]$ is not 
abelian, and so $[SU(n),SU(n)]$ is not abelian.

It remains to prove that $[SU(n),SU(n)]$ is not abelian for
$n= 5,7, 11$.  The argument we now give applies to $SU(p)$ for any prime $p\geq 5$.
Notice that $p$ is regular for $SU(p)$, so it suffices to 
show that $[S^3 \cross S^{2p-3}, SU(p)_{(p)}]$ is nonabelian.
Since $p$ is a regular prime for $SU(p)$, we choose generators 
$\alpha\in \pi_3(SU(p))$, $\beta\in \pi_{2p-3}(SU(p))$ and
$\gamma\in \pi_{2p}(SU(p))$ so that 
$$
\< \alpha', \beta'\> = (p-2)! \gamma' \neq 0 \in \ZZ/p! \otimes \ZZ_{(p)} \cong \ZZ/p.
$$
Therefore  $[S^3\cross S^{2p-3},SU(p)_{(p)}]$ is nonabelian by Lemma \ref{lem:5.3}.

The proof that $[Sp(n),Sp(n)]$ is not abelian for $n\geq 2$ is analogous:
one uses Bott's result for Samelson products in $\pi_*(Sp(n))$ 
\cite[Thm.$\,$2]{Bott}
together with a quasi-regular decomposition for $Sp(n)$.
We omit the details.\halmos
 
\eject
\bcor\hfill
\begin{enumerate}
\item[\rm(a)]  For $n\geq 5$, $Z_\s(SU(n))\neq 0$, and 
$2 \leq t_\s (SU(n)) \leq \left\lceil \log_2(n) \right\rceil$.
\item[\rm(b)]  For $n\geq 2$, $Z_\s(Sp(n))\neq 0$, and  
$2 \leq t_\s(Sp(n)) \leq  \left\lceil 2 \log_2(n+1) \right\rceil $.
\end{enumerate}
\ecor

\pf Proof 
For an H-space $Y$, a commutator in $[X,Y]$ is an element of
$Z_\s(X,Y)$ \cite[Thm.$\,$7]{Strom}.  Thus, if $[Y,Y]$ is nonabelian,
$2\leq t_\s(Y)$.  The upper bound for $t_\s(SU(n))$ comes from
Proposition \ref{prop:Et} since Singhof has shown that
$\cat(SU(n))$\break$ = n$ \cite{Singhof}.  The upper bound on
$t_\s(Sp(n))$ follows from Proposition \ref{prop:killcat}.  \halmos

Schweitzer 
\cite[Ex.$\,$4.4]{Schw}  
has shown that $\cat(Sp(2)) =4$, so
it follows from Proposition \ref{prop:Et} that $t_\s(Sp(2))=2$.

\subsection{Some low dimensional Lie groups}
Here we consider the Lie groups $SU(3)$, $SU(4)$, $SO(3)$ and $SO(4)$
and make estimates of $t_\s$ by either quoting known results or by ad hoc
methods.  We first deal with $SU(3)$ and $SU(4)$.

\bprp\label{prop:su34}
The groups
$[SU(3),SU(3)]$ and $[SU(4),SU(4)]$ are not ab\-elian.
\eprp

\pf Proof
For the group
$\!SU(3)\!$ this follows from results of  Ooshima 
\cite[Thm.$\,$1.2]{Oo}.
For $SU(4)$, observe that the prime $5$ is regular for both $SU(4)$ and $Sp(2)$, so 
$$
SU(4)_{(5)} \equiv (S^3\cross S^5\cross S^7)_{(5)} \quad \mathrm{and}\quad
Sp(2)_{(5)} \equiv (S^3\cross S^7)_{(5)}.
$$
If $[SU(4),SU(4)]$ is abelian, then so is 
$[SU(4)_{(5)},SU(4)_{(5)}] \cong [S^3\cross S^5\cross S^7, SU(4)_{(5)}]$,
and thus $[S^3\cross S^7, SU(4)_{(5)}]$ is abelian.

If $\alpha'\in \pi_3(Sp(2)_{(5)})$ and 
$\beta'\in \pi_7(Sp(2)_{(5)})$ are the images
of generators of $\pi_3(Sp(2))\cong \ZZ$ and  $\pi_7(Sp(2))\cong \ZZ$ then it follows from \cite{Bott} that
$\<\alpha',\beta'\> \neq 0\in \pi_{10}(Sp(2)_{(5)})\cong \ZZ/5!\otimes \ZZ_{(5)} \cong \ZZ/5$.

Now we relate $SU(4)$ to $Sp(2)$ via the fibration
$Sp(2) \maprt{i} SU(4) \maprt{} S^5$.  The exact homotopy sequence of a fibration shows that
$\pi_{10}(i)$ is an isomorphism after localizing at any odd prime. 
Since $i$ is an H-map, 
$$
\< i_*(\alpha'),i_*(\beta')\> = i_*\< \alpha',\beta'\> \neq 0\in \pi_{10}(SU(4)_{(5)}).
$$
Thus $[S^3\cross S^7, SU(4)_{(5)}]$ is not abelian, so $[SU(4),SU(4)]$ cannot be abelian.   \halmos

\bcor\hfill
\begin{enumerate}
\item[\rm(a)]
$Z_\s(SU(3))\neq 0$, and $t_\s(SU(3)) = 2$
\item[\rm(b)]
$Z_\s(SU(4))\neq 0$, and $t_\s(SU(4)) = 2$.
\end{enumerate}
\ecor

\pf Proof
Since the groups
$[SU(n),SU(n)]$ are not abelian for $n=3$ and $4$,
$t_\s(SU(3))$ and $t_\s(SU(4))$ are at least $2$.  But
$\cat(SU(n)) =n$ by \cite{Singhof}, so the reverse inequalities follow from 
Proposition \ref{prop:Et}.\halmos

Next we investigate the nilpotence of $SO(3)$ and $SO(4)$.  This provides us with 
examples of non-simply-connected Lie groups.

\bprp
$Z_\s(SO(3)) = 0$ and $Z_\s(SO(4))\neq 0$.
\eprp

\pf Proof
Since $SO(3)$ is homeomorphic to $\rp^3$, the first assertion follows from
Theorem \ref{thrm:rpn}.  
For the second assertion, recall that $SO(4)$ is homeomorphic to
$S^3\cross SO(3)$.   For notational convenience, we write $X= SO(3)$
and $Y= S^3$.  We show that $Z_\s(X\cross Y)\neq 0$.  Let
$j:X\wdg Y\maprt{} X\cross Y$ be the inclusion and $q:X\cross Y \maprt{} X\smsh Y$
be the quotient map.  Consider 
$$
q^*:[X\smsh Y, X\cross Y]\maprt{} [X\cross Y, X\cross Y].
$$
Notice that $\im(q^*)\sseq Z_\s(X\cross Y)$ because $q\in Z_\s(X\cross Y, X\smsh Y)$, so
$q$ induces a function
$q^{**}: [X\smsh Y, X\cross Y]\maprt{} Z_\s(X\cross Y)$.
Consider the exact sequence of groups
$$
[\s(X\cross Y), X\cross Y]\maprt{\s j^*}[\s(X\wdg Y), X\cross Y]\maprt{}
[X\smsh Y, X\cross Y]\maprt{q^*}[X\cross Y, X\cross Y].
$$
Since $\s j^*$ has a left inverse, $\ker(q^*) =0$.  Thus $q^{**}$ is 
one-one, so it suffices to show that 
$[X\smsh Y , X\cross Y ] \cong [\s^3 SO(3), SO(3)]\oplus [\s^3 SO(3), S^3]$
is nonzero.  This follows from
\cite[Cor.$\,$2.12]{Ya}, 
where it is shown that
$[\s^3 SO(3), S^3]  \cong \ZZ/4 \oplus \ZZ/12$.
\halmos

\bcor
$t_\s(SO(3)) = 1$ and $t_\s(SO(4)) = 2$.
\ecor

\pf Proof
We only have to show that $t_\s(SO(4)) \leq 2$. The remark
following Theorem \ref{thrm:rpn} shows that if $f\in Z_\s(X\cross Y)$
then $f|_{X\wdg Y}=0$, so $q^{**}$ is onto.  Thus $f$ factors 
through a sphere, so we can proceed as in the proof of 
Corollary \ref{cor:thpn}.\halmos

\subsection{The group $\E_\om(Y)$}

We conclude the section by relating $Z_\s(Y)$ to a certain group of homotopy
equivalences of $Y$.  For any space $X$, let $\E_\om(X)\sseq [X,X]$ be the group
of homotopy equivalences $f:X\maprt{} X$ such that 
$\om f = \mathrm{id}$.  This group has been studied by Felix and Murillo
\cite{F-M} and by Pavesic \cite{Pa}.  We note that if $Y$ is an H-space, then the  function
$$
\theta: Z_\s(Y)\maprt{} \E_\om(Y)
$$
defined by $\theta(g) = \mathrm{id} + g$  is a bijection
of pointed sets. In general $\theta$ does not preserve the binary operation
in $Z_\s(Y)$ and $\E_\om(Y)$.  Thus $\E_\om(Y)$ is nontrivial whenever
$Z_\s(Y)$ is nontrivial. 

\bprp 
The groups $\E_\om(Y)$ are nontrivial in the following cases: $Y= SU(n)$, $n\geq 3$;
$Y= Sp(n)$, $n\geq 2$; and $Y= SO(4)$.
The groups $\E_\om(Y)$ are trivial in the following cases: $Y = SU(2)$, $Sp(1)$,
 $SO(2)$ and $SO(3)$.
\eprp

\section{Problems}

In this brief section we list, in no particular order, a number of problems
which extend the previous results or which have been suggested by this material.

\begin{enumerate}
\item[1.] Calculate $t_\F (X)$ for $\F = \S,\M$ or $\s$ and various spaces $X$.
In particular, what is $t_\s (\hp^n)$ for $n>4$, and $t_\s(Y)$ for   compact 
Lie groups $Y$ not considered in Section 6?
\item[2.]
Find general conditions on a space $X$ such that $Z_\s(X,Y) =Z_\S(X,Y)$.  One such
was given in Section 5.  Is $Z_\s(Y)=Z_\S(Y)$ if $Y$ is a compact simply-connected
Lie group without homological torsion, such as $SU(n)$ or $Sp(n)$?
\item[3.]
Find
lower bounds for $t_\F(X)$ in the cases $\F = \S, \M$ or $\s$ in terms of 
other numerical invariants of homotopy type.
\item[4.]
With $\F = \S, \M$ or $\s$, characterize those spaces $X$ such that $Z_\F(X)=0$.
\item[5.]
What is the relation between $\kl_\s(X)$ and $\left\lceil \log_2(\cat(X))\right\rceil$?    
In particular, if $\kl_\s(X)<\infty$, does it follow that $\cat(X)<\infty$?  Notice that
both of these integers are upper bounds for $t_\s(X)$.
\item[6.]
Find an example of a finite H-complex $Y$ such that $Z_\s(Y)\neq Z_\S(Y)$
(see Section 4).  In the notation of Section 6, this would yield a finite 
complex $Y$ for which $\E_\om (Y)\neq \E_\S(Y)$.  Such an example
which is not a finite complex was given in \cite{F-M}.
\item[7.]
Examine $t_\F(X), \kl_\F(X)$ and $\cl_\F(X)$ for various collections $\F$ such
as the collection of $p$-local spheres or the collection of all cell complexes with
at most two positive dimensional cells.
\item[8.]
Investigate the Eckman-Hilton dual of the results of this paper.  One defines 
a map
$f:X\maprt{}Y$ to be \term{$\F$-cotrivial} if 
$f^*=0:[Y,A]\maprt{}[X,A]$
for all $A\in\F$.  One could then study the set $Z^\F(X,Y)$ of all 
$\F$-cotrivial
maps $X\maprt{}Y$ and , in particular, the semigroup $Z^\F(X) = Z^\F(X,X)$.
\end{enumerate}

\Addresses\recd

\end{document}